\documentclass[reqno]{amsart}
\usepackage[T1]{fontenc}
\usepackage{amsmath,amsthm,amssymb,enumitem,times,xcolor,bm,esint,mathtools,dsfont}
\usepackage{leftidx}
\usepackage[mathscr]{eucal}
\usepackage{harpoon}
\usepackage{todonotes}
\usepackage{tikz}
\usepackage[titletoc]{appendix} 
\usepackage[
linktocpage=true,colorlinks,citecolor=blue,linkcolor=blue,urlcolor=blue]{hyperref}
\allowdisplaybreaks[4]
\providecommand{\U}[1]{\protect \rule{.1in}{.1in}}

\hypersetup{citecolor=red,linkcolor=red}
\usepackage{geometry}
\newtheorem{theorem}{Theorem}[section] 
\newtheorem{corollary}[theorem]{Corollary} 
\newtheorem{definition}[theorem]{Definition} 
\newtheorem{lemma}[theorem]{Lemma} 
\newtheorem{proposition}[theorem]{Proposition}
\newtheorem{theoremA}{Theorem} 
\theoremstyle{remark}
 
\numberwithin{equation}{section}  

\def\be{\begin{equation}}
\def\ee{\end{equation}}
\def\lab{\label}
\def\AR{{\mathcal A}} 
\def\BR{{\mathcal B}}
\def\CR{{\mathcal C}}
\def\DR{{\mathcal D}}
\def\ER{{\mathcal E}}

\def\HR{{\mathcal H}}

\def\LR{{\mathcal L}}
\def\MR{{\mathcal M}}

\def\OR{{\mathcal O}}

\def\QR{{\mathcal Q}}
\def\SR{{\mathcal S}}
\def\TR{{\mathcal T}}

\def\R{\mathbb R}
\def\N{\mathbb N}

\def\al{\alpha}

\newcommand{\ve}{\varepsilon}
\newcommand{\vs}{\varsigma}

\def\iy{\infty}
\def\la{\lambda}

\def\vp{\varphi}
\def\pa{\partial}
\def\ti{\tilde}
\def\wt{\widetilde}
\def\ov{\overline}
\def\f{\frac}

\def\sbr#1{\left(#1\right)} 
\def\lbr#1{\left\{#1\right\}} 
\def\abr#1{\left\langle#1\right\rangle} 
\def\abs#1{\left\lvert#1\right\rvert} 
\def\nm#1{\left\|#1\right\|} 
\def\vs{\vskip .2in} 
\def\rd{\mathrm d}   
\def\loc{{\operatorname{\rm loc}}} 
\def\id{\operatorname{\rm id}} 
\def\bxi{{\boldsymbol{\xi}}}
\def\bX{{\boldsymbol{X}}}
\def\bmu{{\boldsymbol{\mu}}}
\def\bal{{\boldsymbol{\alpha}}}

\begin{document}
\theoremstyle{plain}

\title{Symmetry and uniqueness of low-energy nodal solutions for the Lane--Emden equation in the disk}

\author{Zhijie Chen}\address{Department of Mathematical Sciences, Yau Mathematical Sciences Center, Tsinghua University, Beijing, China}\email{zjchen2016@tsinghua.edu.cn}
\author{Houwang Li}\address{School of Mathematics and Statistics, Beijing Institute of Technology, Beijing 100081,
P.R. China}\email{houwangli@bit.edu.cn}
\author{Hanqing Zhao}
\address{School of Mathematics and Statistics, Wuhan University, Wuhan 430072, China}
\email{zhq25math@whu.edu.cn}

\begin{abstract}
We study nodal solutions of the Lane--Emden equation in the disk
\[
-\Delta u=|u|^{p-1}u\quad\text{in }B_1(0),\qquad
u=0\quad\text{on }\partial B_1(0).
\]
Esposito, Musso and Pistoia (Proc.LMS, 2007) constructed two-peak nodal solutions for large $p$ and conjectured that they minimize the energy among all nodal solutions.
We confirm their conjecture in the disk by proving that any two low-energy nodal families agree, up to rotation, for all sufficiently large $p$.
As $p\to\infty$, the ordered pair of peak locations approaches a rotational orbit of degenerate critical points of the associated Kirchhoff--Routh function.
The rotational degeneracy prevents a direct application of existing local uniqueness results.
We first establish axial symmetry by a rotating-plane argument and then compute the Morse index and the Leray--Schauder degree in a symmetric subspace excluding the rotational mode.
Each low-energy nodal solution is also antisymmetric about the diameter perpendicular to its symmetry axis. 
In the full space, its Morse index is two and its linearized kernel is spanned by the rotational derivative.
\end{abstract}

\maketitle

\section{Introduction}
\subsection{Motivation}
In a smooth bounded domain $\Omega\subset\R^2$, the Lane--Emden problem is
\begin{equation}\label{equ-1****}\begin{cases}
    -\Delta u = |u|^{p-1}u, & \text{in } \Omega, \\
    u = 0, & \text{on } \partial \Omega,
\end{cases}\end{equation}
where $1<p<\iy$. Families of solutions to \eqref{equ-1****} may concentrate as $p\to\infty$.
The concentration of positive solutions with bounded rescaled energy has been studied extensively. See, for example, \cite{Adimurthi-Grossi-2004,Ren-Wei-94,Ren-Wei-96,DeMarchis-Ianni-Pacella-2017-2,DeMarchis-Grossi-Ianni-Pacella-2018,Thizy-2019,Kamburov-Sirakov-2018,DeMarchis-Grossi-Ianni-Pacella-2019,Grossi-Ianni-Luo-Yan-2022}.
The concentration of low-energy nodal solutions is described in \cite{Grossi-Grumiau-Pacella-2013,Chen-Cheng-Zhao-2024}. Other nodal concentration patterns are studied in \cite{DeMarchis-Ianni-Pacella-2015,DeMarchis-Ianni-Pacella-2016,DeMarchis-Ianni-Pacella-2017,Gladiali-Stegel-2022,Grossi-Grumiau-Pacella-2014,Ianni-Saldana-2021}.

A family of nodal solutions $\{u_p\}_{p>1}$ of \eqref{equ-1****} is called a family of {\it low-energy nodal solutions} if
\be\lab{low-energy} \lim_{p\to+\iy} p\int_\Omega |u_p|^{p+1} \rd x=8\pi e\times2. \ee
A {\it least-energy nodal solution} minimizes the energy among all nodal solutions. 
By Grossi, Grumiau and Pacella \cite{Grossi-Grumiau-Pacella-2013}, every family of least-energy nodal solutions satisfies \eqref{low-energy}.
Thus the low-energy class includes least-energy nodal families. 
To describe the location of positive and negative peaks of a low-energy nodal family, we need the Kirchhoff--Routh type function $\Psi:\MR=\Omega^2\setminus\{(x,y):x=y\}\to\R$, defined by
\be\lab{Psi} \Psi(a_1,a_2):=\sum_{i=1}^2\Psi_i(a_1,a_2)=-R(a_1)-R(a_2)+2G(a_1,a_2), \ee
where $G$ is the Dirichlet Green function of $-\Delta$ in $\Omega$. 
We write
    \[G(x,y)=\f{1}{2\pi}\log\f{1}{|x-y|}+H(x,y),\]
and set the Robin function $R(x)=H(x,x)$. Note that the sign convention for the Robin function varies in the literature. 
For later use, we also set
\be\lab{Psi_i} \Psi_i(a_1,a_2)=-R(a_i)+G(a_1,a_2). \ee
Along a subsequence, the global maximum and minimum points converge to distinct points $x^+,x^-\in\Omega$ satisfying
    \[\nabla\Psi(x^+,x^-)=0.\]
Grossi, Grumiau and Pacella \cite{Grossi-Grumiau-Pacella-2013} obtained this description under additional assumptions, including a comparison condition on the two peak heights.
Chen, Cheng and Zhao \cite{Chen-Cheng-Zhao-2023} gave an argument for this concentration description without the peak-height comparison. They subsequently proved that the comparison holds for every low-energy nodal family \cite{Chen-Cheng-Zhao-2024}.

\vskip0.05in
The local uniqueness and nondegeneracy of concentrating solutions require more than the location of the limiting peaks. 
For positive one-peak solutions, the Morse index, nondegeneracy and local uniqueness were obtained in \cite{DeMarchis-Grossi-Ianni-Pacella-2019} when the peak converges to a nondegenerate critical point of the Robin function.
For positive multi-peak solutions, the local uniqueness results in \cite{Grossi-Ianni-Luo-Yan-2022,Ianni-Luo-Yan-2025} also require nondegeneracy of the corresponding concentration configuration.
For low-energy nodal solutions, Chen, Cheng and Zhao \cite{Chen-Cheng-Zhao-2024} proved nondegeneracy when $(x^+,x^-)$ is a nondegenerate critical point of $\Psi$.

In the disk, Chen, Cheng and Zhao \cite{Chen-Cheng-Zhao-2024} observed rotational degeneracy for least-energy nodal solutions and their limiting configuration.
As shown in Section~\ref{sec-2}, all critical point configurations of $\Psi$ are rotations of a single antipodal pair. 
At each configuration, the Hessian has a one-dimensional kernel tangent to the rotational orbit and is positive definite on its orthogonal complement, see \eqref{Hessian}. 
Thus the concentration configuration is degenerate, and a nonradial solution also has a rotational mode in its linearized kernel. 
The uniqueness and nondegeneracy statements must therefore be formulated modulo rotations.

The question is whether low-energy nodal solutions in the disk are unique modulo rotations for large $p$ and whether the rotational mode spans the linearized kernel.
The same low-energy class also arises in the construction of Esposito, Musso and Pistoia \cite{Esposito-Musso-Pistoia-2007}.
Their existence result uses a minimum of $\Psi$ and does not require the minimum to be nondegenerate.

\begin{theoremA}[\cite{Esposito-Musso-Pistoia-2007}]\label{Esposito-Musso-Pistoia}
There exists $p_0 > 0$ such that, for any $p\geq p_0$, problem \eqref{equ-1****} has a pair of sign-changing solutions $\pm u_p$ such that
    \[p|u_p|^{p-1}u_p \mathrm{d}x\rightharpoonup 8\pi\sqrt{e} \left(\delta_{\xi_1^*}-\delta_{\xi_2^*}\right), \]
weakly as signed measures in $\Omega$ as $p\to+\infty$ for some $\xi_1^*,\xi_2^*\in\Omega$ with $\xi_1^*\neq\xi_2^*$, such that $\Psi(\xi_1^*,\xi_2^*)=\min_{\MR}\Psi$.  These solutions satisfy
    \[\lim_{p\to+\infty}p\int_\Omega |u_p|^{p+1}\,\rd x=16\pi e.\]
\end{theoremA}

Esposito, Musso and Pistoia conjectured that, for sufficiently large $p$, the solutions constructed in Theorem~\ref{Esposito-Musso-Pistoia} minimize the energy among all nodal solutions of \eqref{equ-1****}.
In the disk, the classification below identifies these constructed solutions with the least-energy nodal solution, up to rotations, for sufficiently large $p$, and hence answers their conjecture in the disk.

Bartsch and Weth \cite{Bartsch-Weth-2003} proved that least-energy nodal solutions have Morse index two and exactly two nodal domains. On radial domains, Bartsch, Weth and Willem \cite{Bartsch-Weth-Willem-2005} established foliated Schwarz symmetry for least-energy nodal solutions. For the disk and $p\ge2$, the result of Pacella and Weth \cite{Pacella-Weth-2007} also gives this symmetry for solutions of Morse index at most two. Aftalion and Pacella \cite{Aftalion-Pacella-2004} showed that least-energy nodal solutions on a ball are nonradial and that the closure of their nodal set meets the boundary.

In the disk, Grossi, Grumiau and Pacella \cite{Grossi-Grumiau-Pacella-2014} proved antisymmetry of the limiting profile of $pu_p$ under the assumptions of Morse index two and comparable peak heights. They also discussed whether least-energy nodal solutions themselves are antisymmetric. Fazekas, Pacella and Plum \cite{Fazekas-Pacella-Plum-2022} formulated this question in terms of odd reflection across a diameter. Exact antisymmetry on a ball when $p$ is close to $1$ was proved in \cite{Bonheure-Bouchez-Grumiau-VanSchaftingen-2008}, and uniqueness up to rotations and sign in that regime was proved in \cite{Grumiau-Troestler-2010}. Here we show that, for every sufficiently large $p$, the entire low-energy class in the disk forms a single rotational orbit and the interior nodal set of each solution is exactly one diameter.

\subsection{Setting and Main Results}
We now take $\Omega$ to be the unit disk and consider the Lane--Emden problem:
\begin{equation}\label{equ-1}\begin{cases}
    -\Delta u = |u|^{p-1}u, & \text{in } B_1(0), \\
    u = 0, & \text{on } \partial B_1(0),
\end{cases}\end{equation}
where $B_1(0) = \{x \in \mathbb{R}^2 : |x| < 1\}$ is the unit disk in $\mathbb{R}^2$ and $p > 1$.
Let $J_p:H_0^1(B_1(0))\to\R$ be the associated energy functional
\be\lab{J_p} J_p(u)=\f12\int_{B_1(0)} |\nabla u|^2\rd x-\f{1}{p+1}\int_{B_1(0)} |u|^{p+1} \rd x. \ee
Its critical set is
\be\lab{CR} \CR_p=\{u\in H_0^1(B_1(0)):~J_p'(u)=0\}. \ee
Thus $\CR_p$ is the set of solutions of \eqref{equ-1}. 
In order to quantitatively characterize low-energy solutions, we define the open set
\be\lab{define-e_p}     \ER_p=\left\{u\in H_0^1(B_1(0)):0<16\pi e-p\int_{B_1(0)}|u|^{p+1}\,\rd x <48\pi e\,\frac{\log p}{p}\right\}. \ee
The positive solution $v_p$ of \eqref{equ-1} is unique by \cite{Gidas-Ni-Nirenberg-1979}. It is the positive least-energy solution, so $p\int_{B_1(0)}|v_p|^{p+1}\,\rd x=8\pi e+o(1)$ by the standard least-energy asymptotics \cite{Ren-Wei-94,Adimurthi-Grossi-2004}. 
Consequently, every solution in $\overline{\ER_p}$ is nodal for sufficiently large $p$.
Moreover, it is easy to see that every sequence of such solutions with $p\to\infty$ satisfies the energy limit in \eqref{low-energy}. 
Conversely, Lemma~\ref{lem-mass-energy} places every low-energy nodal family in $\ER_p$ for all sufficiently large $p$.

Our first result determines the extrema and the reflection symmetries of these solutions.

\begin{theorem}[Symmetry]\lab{thm-symmetry}
There exists $p_0>1$ such that, for every $p>p_0$ and every $u_p\in\ER_p\cap\CR_p$, $u_p$ has a unique global maximum point $x_{p,1}$ and a unique global minimum point $x_{p,2}$. Moreover, the points $x_{p,1},x_{p,2}$ and the origin are collinear. The solution $u_p$ is symmetric about this line and antisymmetric about the perpendicular line through the origin.
Its nodal set in $B_1(0)$ is precisely the diameter on that perpendicular line.
\end{theorem}

The same solution class consists of a single rotational orbit, and the linearized kernel at each solution is spanned by the rotational mode.

\begin{theorem}[Uniqueness]\lab{thm-uniqueness}
There exists $p_0>1$ such that, for every $p>p_0$,
    \[\sharp\bigl((\ER_p\cap\CR_p)/\!\sim\bigr)=1,\]
where $\sim$ denotes equivalence under rotations about the origin. Thus, up to rotation, the solution in $\ER_p\cap\CR_p$ is unique for all $p>p_0$.
Moreover, each solution $u_p\in\ER_p\cap\CR_p$ has Morse index $2$ in $H_0^1(B_1(0))$, and it is nondegenerate modulo rotations in the sense that the linear problem 
\be\lab{linear problem}\begin{cases}
    -\Delta \phi=p|u_p|^{p-1}\phi,\quad\text{in}~B_1(0),\\
    \phi=0,\quad\text{on}~\pa B_1(0),
\end{cases}\ee
has kernel $\R \pa_\theta u_p$, where $\pa_\theta u_p=-x_2\pa_{x_1}u_p+x_1\pa_{x_2}u_p$. 
\end{theorem}

We explain the strategy of the proof.
Wei and Winter \cite{Wei-Winter-2005} established symmetry and uniqueness for two-peak nodal solutions of a singularly perturbed equation on a ball. Following their rotating-plane argument, we establish axial symmetry, so that every solution in $\ER_p\cap\CR_p$ can be rotated into the subspace $\HR_e\subset H_0^1(B_1(0))$ of functions even in $x_2$.
Then we use the spectral approach developed in \cite{DeMarchis-Grossi-Ianni-Pacella-2019,Ianni-Luo-Yan-2025} to study the eigenvalue problem 
\[\begin{cases}
    -\Delta v=\la p|u_p|^{p-1}v \quad\text{in}~B_1(0),\\
    v=0 \quad\text{on}~\pa B_1(0),\\
    v\in\HR_e.
\end{cases}\]
The Morse index in $\HR_e$ counts the eigenvalues below $1$, while the eigenvalue $1$ corresponds to the linearized kernel. 
Applying the min--max principle to cutoffs of $u_p$ supported near the two peaks gives $\lambda_{p,1},\lambda_{p,2}\to0$, and adding cutoffs of $\partial_{x_1}u_p$ to the test space gives $1+o(1)=\lambda_{p,3}\le\lambda_{p,4}\leq 1+O(\varepsilon_p^2)$. 
Then, by the Green function expansions and local Pohozaev identities, we express the leading shifts of $\lambda_{p,3}$ and $\lambda_{p,4}$ in terms of the Hessian of $\Psi$ restricted to displacements along the axis.
This restricted Hessian is positive definite, so both eigenvalues lie strictly above $1$ for large $p$. 
It follows that each solution in $\ER_p\cap\CR_p\cap\HR_e$ has Morse index $2$ and a trivial linearized kernel in $\HR_e$.
These solutions are isolated, and each contributes local Leray--Schauder degree $(-1)^2=+1$.
Thus the total degree counts the solutions in $\ER_p\cap\CR_p\cap\HR_e$.

We then compute the total degree by a finite-dimensional reduction, which follows the framework of \cite{Esposito-Musso-Pistoia-2007,Ianni-Luo-Yan-2025}.
We write \eqref{equ-1} as $(\id-\TR_p)(u)=0$, where $\TR_p:H_0^1(B_1(0))\to H_0^1(B_1(0))$ is defined by
\be\lab{TR} \TR_p(u)=(-\Delta)^{-1}(|u|^{p-1}u). \ee
For fixed $p$, the operator $\TR_p$ is compact on $H_0^1(B_1(0))$. 
To compute the degree, we work in the bounded open set
\be\lab{define-d_p} \DR_p=\ER_p\cap\{u\in H_0^1(B_1):\|u\|<r_p\}\quad\text{with}\quad r_p=2\sqrt{\frac{16\pi e}{p}}.\ee
It is easy to check that $\ER_p\cap\CR_p=\DR_p\cap\CR_p$ and
    \[(I-\TR_p)(u)\ne0\qquad\text{for every }u\in\partial\DR_p,\]
so that the degree of $I-\TR_p$ is well defined in $\DR_p\cap\HR_e$.
By representing each solution $u\in\ER_p\cap\CR_p\cap\HR_e$ as a signed two-peak profile with two amplitude parameters $\boldsymbol\al$, two axial location parameters $\boldsymbol\xi$ and a small remainder $w$, we can reduce the computation of the Leray--Schauder degree to that of an operator $\AR_p$.
By homotopy invariance, we reduce the degree calculation for $\AR_p$ to the amplitude equations, the axial gradient of $\Psi$ and the identity on the remainder space. The first two blocks form a finite-dimensional problem.
The amplitude block and the identity on the remainder space each have degree $+1$, and the degree of the axial gradient of $\Psi$ is again determined by the restricted Hessian. 
Thus, the total degree is $2$, which gives exactly two solutions in $\ER_p\cap\CR_p\cap\HR_e$, i.e., $u_p$ and $-u_p$.
Rotation through $\pi$ exchanges the two solutions, proving uniqueness up to rotations. 

\vskip0.1in
Least-energy nodal families and the family in Theorem~\ref{Esposito-Musso-Pistoia} are low-energy families. Theorem~\ref{thm-uniqueness} therefore identifies them in the disk and gives the following consequences. In particular, it answers the conjecture of Esposito, Musso and Pistoia affirmatively for the disk.

\begin{corollary}[Least-energy nodal solution]\lab{cor-1}
There exists $p_0>1$ such that, for every $p>p_0$, equation \eqref{equ-1} has a unique least-energy nodal solution up to rotations.
\end{corollary}

\begin{corollary}[The conjecture in the disk]\lab{cor-2}
Let $\{u_p\}_{p>1}$ be a family of low-energy nodal solutions of \eqref{equ-1}. There exists $p_0>1$ such that, for every $p>p_0$, $u_p$ coincides, up to a rotation, with the least-energy nodal solution.
\end{corollary}

\noindent{\bf Outline.}
Section~\ref{sec-2} describes the concentration phenomenon as $p\to+\iy$ and records some asymptotic estimates. 
Section~\ref{sec-3} proves axial symmetry and angular monotonicity. 
Sections~\ref{sec-4} and~\ref{sec-degree} compute the Morse index and the degree in the even subspace. 
Section~\ref{sec-6} uses these computations to prove antisymmetry, uniqueness and the full-space nondegeneracy statement.

\vskip0.1in
\noindent{\bf Notation.}
Throughout the paper, $C$ denotes a positive constant that is independent of $p$ and may vary from line to line. 
We write $B_r(x)=\{y\in\R^2:|y-x|<r\}$ and sometimes write $B_r=B_r(0)$. 
We use $\|\cdot\|$ for the Dirichlet norm on $H_0^1(B_1)$.
We set $u^+=\max\{u,0\}$ and $u^-=\max\{-u,0\}$.
For a function of two variables, $\partial_i$ and $D_i$ denote differentiation in the first and second variables, respectively.
Limits of concentration points and rescaled functions are understood after passage to a subsequence unless stated otherwise.

\subsection*{AI Declaration}
The authors developed the main ideas of the proof independently of AI tools.
They used OpenAI's ChatGPT and Codex to edit the English and the LaTeX writing, format references, and review mathematical notation, calculations and proofs.
These tools also assisted in deriving and presenting selected proof steps, including the fixed-space formulation and the local-degree comparison in Section~\ref{sec-degree}.
The authors take full responsibility for the mathematical arguments, citations and final text.

\vskip0.2in
\section{Basic results}\lab{sec-2}
\subsection{Critical points of the Kirchhoff--Routh type function}
In the unit disk, the Green function is
\begin{equation}\label{eq-green}
    G(z_1,z_2)=\f{1}{2\pi}\log\f{1}{|z_1-z_2|}+\f{1}{2\pi}\log|z_1\bar z_2-1|,
\end{equation}
where $z_1,z_2\in B_1(0)$ are complex numbers and $\bar z$ denotes the conjugate of $z$. Then the Robin function is
    \[R(z)=\f{1}{2\pi}\log(1-|z|^2)\]
and the Kirchhoff--Routh type function is
\be\lab{Psi-1} \Psi(z_1,z_2)=-\f{1}{2\pi}\log\sbr{(1-|z_1|^2)(1-|z_2|^2)}-\f{1}{2\pi}\log|z_1-z_2|^2+\f{1}{2\pi}\log|z_1\bar z_2-1|^2, \ee
with $z_1,z_2\in B_1(0)$ and $z_1\neq z_2$. 

\begin{lemma}
The set of critical points of $\Psi$ is
    \[\lbr{(\sqrt{\sqrt5-2}e^{\sqrt{-1}\theta},-\sqrt{\sqrt5-2}e^{\sqrt{-1}\theta}):\theta\in[0,2\pi)}.\]
\end{lemma}
\begin{proof}
Differentiating $\Psi$ gives the critical point equations 
\be\begin{cases}
    2\pi\pa_{z_1}\Psi(z_1,z_2)=-\f{\bar z_1}{|z_1|^2-1}-\f{1}{z_1-z_2}+\f{\bar z_2}{z_1\bar z_2-1}=0,\\
    2\pi\pa_{z_2}\Psi(z_1,z_2)=-\f{\bar z_2}{|z_2|^2-1}+\f{1}{z_1-z_2}+\f{\bar z_1}{\bar z_1z_2-1}=0.
\end{cases}\ee
These equations imply $z_1\neq0$ and $z_2\neq0$.
Set $z_1=se^{\sqrt{-1}w}$ and $z_2=re^{\sqrt{-1}(w+\vp)}$ for $s,r\in(0,1)$ and $w,\vp\in[0,2\pi)$. 
For distinct points, the rotational invariance of $\Psi$ gives
\[\begin{aligned} 
    F(s,r,\vp)&=-2\pi\Psi(se^{\sqrt{-1}w},re^{\sqrt{-1}(w+\vp)})\\
        &=\log\sbr{(1-s^2)(1-r^2)}+\log(s^2+r^2-2sr\cos\vp)-\log(s^2r^2+1-2sr\cos\vp).
\end{aligned}\]
Thus $(s,r,\vp)$ is a critical point of $F$ if and only if $(se^{\sqrt{-1}w},re^{\sqrt{-1}(w+\vp)})$ is a critical point of $\Psi$ for every $w\in[0,2\pi)$.
Differentiation shows that $\pa_\vp F(s,r,\vp)=0$ if and only if $\vp\in\{0,\pi\}$. We set $T=(0,1)\times(0,1)$ and prove the following claim. $F(s,r,0):T\setminus\{s=r\}\to \R$ has no critical point, and $F(s,r,\pi):T\to \R$ has a unique critical point $s=r=\sqrt{\sqrt5-2}$.

If $\vp=0$, we change variables and set $\xi=(s-r)^2$, $\eta=sr$. Notice that 
    \[\det\sbr{\f{\pa(\xi,\eta)}{\pa(s,r)}}=2(s^2-r^2)=0,\quad\text{if and only if $s=r$ in $T$}.\]
An explicit evaluation shows that 
    \[\pa_\eta F(s(\xi,\eta),r(\xi,\eta),0)=\f{-2\xi}{(1-\eta)((1-\eta)^2-\xi)}\neq0,\]
so that any critical point of $F(s,r,0)$ must satisfy $s=r$, and hence $F(s,r,0)$ has no critical point in $T\setminus\{s=r\}$.

If $\vp=\pi$, we change variables and set $\xi'=(s+r)^2$, $\eta'=sr$. Notice that 
    \[\det\sbr{\f{\pa(\xi',\eta')}{\pa(s,r)}}=2(s^2-r^2)=0,\quad\text{if and only if $s=r$ in $T$}.\]
An explicit evaluation shows that 
    \[\pa_{\eta'} F(s(\xi',\eta'),r(\xi',\eta'),\pi)=\f{2\xi'}{(1+\eta')((1+\eta')^2-\xi')}\neq0,\]
so that any critical point of $F(s,r,\pi)$ must satisfy $s=r$. At this point, the equations give $\pa_sF(t,t,\pi)=\pa_rF(t,t,\pi)=0$ if and only if $t=\sqrt{\sqrt5-2}$.
\end{proof}

After a rotation, we may place the critical pair on the $x_1$-axis.
\be\lab{X}     X_1=(\sqrt{\sqrt5-2},0), \quad X_2=(-\sqrt{\sqrt5-2},0). \ee
Direct differentiation gives the Hessian of $\Psi$ at $(X_1,X_2)$,
\be\lab{Hessian}
    H_\Psi(X_1,X_2) =\f{1}{8\pi}\sbr{ \begin{matrix}
    13+5\sqrt5 & 0 & -7-3\sqrt5 & 0\\
    0 & 1+\sqrt5 & 0 & 1+\sqrt5\\
    -7-3\sqrt5 & 0 & 13+5\sqrt5 & 0\\
    0 & 1+\sqrt5 & 0 & 1+\sqrt5
    \end{matrix} },
\ee
and the determinant equals 0. This shows that $(X_1,X_2)$ is a degenerate critical point of $\Psi$.

\subsection{Asymptotic estimates for low-energy nodal solutions}\label{sec2.2}
Let $u_p$ be a family of low-energy nodal solutions of \eqref{equ-1}. 
We recall in (P1)--(P2) the planar concentration results in \cite{Grossi-Grumiau-Pacella-2013,Chen-Cheng-Zhao-2023}, as summarized in \cite[Section~2]{Chen-Cheng-Zhao-2024}.
After passing to a subsequence, the following properties hold. All estimates are stated for a fixed family; uniform conclusions for $\ER_p$ follow by the sequential contradiction argument used below.
\begin{itemize}
\item[(P1)] Choose a global maximum point $x_{p,1}$ and a global minimum point $x_{p,2}$ of $u_p$, so that
\[u_p(x_{p,1})=\max_{B_1(0)}u_p=\nm{u_p^+}_{L^\iy(B_1(0))},\quad u_p(x_{p,2})=\min_{B_1(0)}u_p=-\nm{u_p^-}_{L^\iy(B_1(0))}.\]
After passage to a subsequence, write
\[x_{\iy,i}=\lim_{p\to+\iy}x_{p,i},\quad\forall~i=1,2.\]
The pair $(x_{\iy,1},x_{\iy,2})$ is a critical point of the Kirchhoff--Routh-type function $\Psi$ defined in \eqref{Psi}, i.e.,
    \[\nabla\Psi(x_{\iy,1},x_{\iy,2})=0.\]

\item[(P2)] For $i=1,2$, we define
\be\lab{ve_pi} \ve_{p,i}=\sbr{p|u_{p}(x_{p,i})|^{p-1}}^{-\f12}\to0, \ee
and
\be\lab{w_pi} w_{p,i}(x)=p\sbr{\f{u_p(x_{p,i}+\ve_{p,i}x)}{u_p(x_{p,i})}-1}, \quad x\in\Omega_{p,i}=\f{B_1(0)-x_{p,i}}{\ve_{p,i}}.\ee
We have $w_{p,i}\to U$ in $C_\loc^2(\R^2)$, where 
\be\lab{U} U(x)=-2\log\sbr{1+\f{1}{8}|x|^2} \ee 
is the radial solution of $-\Delta U=e^U$ in $\R^2$ with $U(0)=0$ and $\int_{\R^2}e^U\rd x=8\pi$.

\item[(P3)] By \cite[Lemmas~3.1, 3.2 and 3.6 and Corollary~3.4]{Chen-Cheng-Zhao-2024}, there exist a small $\rho_0>0$ and a large $p_0>1$ such that, for every $p>p_0$,
\be\lab{3-6-2} u_p\ge\f1p\quad \text{in}~B_{4\rho_0}(x_{p,1}),\quad u_p\le-\f1p \quad \text{in}~B_{4\rho_0}(x_{p,2}), \ee
and
\be 0\le 1+\f{w_{p,i}}{p}\le 1 \quad \text{in}~B_{\f{4\rho_0}{\ve_{p,i}}}(0),\quad\forall~i=1,2. \ee
Moreover, for any $i=1,2$ and any small constant $\delta>0$, there exists $C_\delta>0$ such that
\be\lab{est-wpi-1} w_{p,i}(x)\le (4-\delta)\log\f{1}{|x|}+C_\delta \quad \text{in}~B_{\f{2\rho_0}{\ve_{p,i}}}(0), \ee
\be\lab{est-wpi-2} |w_{p,i}(x)|\le C\log(2+|x|) \quad \text{in}~B_{\f{2\rho_0}{\ve_{p,i}}}(0), \ee
and
\be\lab{est-wpi-3} \abs{1+\f{w_{p,i}(x)}{p}}^p\le \f{C_\delta}{1+|x|^{4-\delta}} \quad \text{in}~B_{\f{2\rho_0}{\ve_{p,i}}}(0). \ee

\item[(P4)] Write $\mathbf x_p=(x_{p,1},x_{p,2})$ and $\mathbf x_\infty=(x_{\infty,1},x_{\infty,2})$.
For $i=1,2$ and every sufficiently small $\delta>0$, \cite[Theorem~1.3 and (3.54)]{Chen-Cheng-Zhao-2024} give
\be
\ve_{p,i}=e^{-p/4}\left(
 e^{-2\pi\Psi_i(\mathbf x_p)-\frac32\log2-\frac34}
 +O(p^{-1+\delta})\right),
\ee
and
\be\lab{3-12-11}
u_p(x_{p,i})=(-1)^{i-1}\sqrt e\left(
1-\frac{\log p}{p}
+\frac{4\pi\Psi_i(\mathbf x_p)+3\log2+2}{p}\right)
+O(p^{-2+\delta}).
\ee
Here $\Psi_i$ is defined in \eqref{Psi_i}.
By taking 
\be\lab{ve_p} \ve_p=\ve_{p,1}e^{2\pi\Psi_1(\mathbf{x}_\iy)}, \ee 
we obtain
\be\lab{ve_p-1} \lim_{p\to+\iy}\f{\ve_{p,2}e^{2\pi\Psi_2(\mathbf{x}_\iy)}}{\ve_p}=1.\ee 

\item[(P5)] For any fixed $\rho\in(0,\rho_0)$, the exterior expansion and signed masses satisfy the following estimates (see \cite[Lemma~3.5 and Remark~3.10]{Chen-Cheng-Zhao-2024}).
\be\lab{out-u}  u_p(x)=\sum_{i=1}^2C_{p,i}G(x_{p,i},x)+o\sbr{\f{\ve_p}{p}}  \quad\text{in}~C_\loc^1\big(\ov{B_1(0)}\setminus \cup_{i=1}^2B_{2\rho}(x_{p,i})\big), \ee 
where 
\be\lab{Cpi} C_{p,i}=\int_{B_{\rho}(x_{p,i})} |u_p|^{p-1}u_p \rd x=(-1)^{i-1}\f{8\pi\sqrt e}{p}\sbr{1+O\sbr{\f{\log p}{p}}}. \ee

\item[(P6)] Fix $\rho\in(0,\rho_0)$. For $i=1,2$ and any sufficiently small $\delta>0$, \cite[Proposition~3.9]{Chen-Cheng-Zhao-2024} gives
\be\lab{3-12-10} \int_{B_{\rho/\ve_{p,i}}(0)} \abs{1+\f{w_{p,i}(x)}{p}}^p\rd x=8\pi-\f{24\pi}{p}+O\sbr{\f{1}{p^{2-\delta}}}, \ee
and
\be \int_{B_{\rho/\ve_{p,i}}(0)} \abs{1+\f{w_{p,i}(x)}{p}}^{p+1}\rd x=8\pi-\f{40\pi}{p}+O\sbr{\f{1}{p^{2-\delta}}}. \ee
\end{itemize}

The first correction to the rescaled profile is the radial solution $w_0$ of
\be\lab{3-12-0}
-\Delta w_0-e^Uw_0=-\frac12U^2e^U\quad\text{in }\R^2,
\qquad w_0(0)=0.
\ee
The next correction is needed for the finite-dimensional approximation.
\begin{lemma}\label{lem-profile-corrections}
Let $\rho\in(0,\rho_0)$ and $\tau\in(0,1)$ be fixed. For $i=1,2$,
\be\lab{est-wpi-4}
\left|w_{p,i}(y)-U(y)-\frac{w_0(y)}p-\frac{w_1(y)}{p^2}\right|
\le\frac{C(1+|y|^\tau)}{p^3},
\qquad |y|\le\frac\rho{\ve_{p,i}},
\ee
where $w_1$ is radial, $w_1(0)=0$, and
\be\label{eq-w1-correct}
-\Delta w_1-e^Uw_1
=e^U\left(\frac12w_0^2-\frac12U^2w_0-Uw_0
+\frac13U^3+\frac18U^4\right)\quad\text{in }\R^2.
\ee
On each fixed compact subset of $\R^2$, the remainder in \eqref{est-wpi-4} is also $O(p^{-3})$ in $C^2$.
\end{lemma}
\begin{proof}
The first two weighted remainders are estimated in \cite[Propositions~3.7 and 3.8]{Chen-Cheng-Zhao-2024}. For the third weighted remainder, we follow \cite[Proposition~5.1 and Lemmas~5.6--5.8]{Ianni-Luo-Yan-2025} and the local argument of \cite[Proposition~3.5]{Grossi-Ianni-Luo-Yan-2022}, verifying the localization for the two signed peaks. By (P3),
\[
 -\Delta w_{p,i}=\left(1+\frac{w_{p,i}}p\right)^p,
 \qquad 1+\frac{w_{p,i}}p>0
 \quad\hbox{on }B_{2\rho/\ve_{p,i}}(0),
\]
for both $i=1,2$. Set
\[\begin{aligned}
v_{p,i}=p(w_{p,i}-U),\qquad
k_{p,i}=p(v_{p,i}-w_0),\qquad
q_{p,i}=p(k_{p,i}-w_1),
\end{aligned}\]
then we have
\[
v_{p,i}(0)=k_{p,i}(0)=q_{p,i}(0)=0,\qquad
\nabla v_{p,i}(0)=\nabla k_{p,i}(0)=\nabla q_{p,i}(0)=0.
\]
The source outside the nodal ball contributes the harmonic function
\[
h_{p,i}(x)=\int_{B_1\setminus B_{2\rho}(x_{p,i})}
G(x,z)|u_p(z)|^{p-1}u_p(z)\,\rd z.
\]
By H\"older's inequality, for each fixed $k$ we obtain
\[\|h_{p,i}\|_{C^k(B_\rho(x_{p,i}))}\le\frac{C_k}{p}.\]
For $M_{p,i}=u_p(x_{p,i})$, write
\[\begin{aligned}
H_{p,i}(y)=\frac p{M_{p,i}}
\bigl(h_{p,i}(x_{p,i}+\ve_{p,i}y)-h_{p,i}(x_{p,i})\bigr).
\end{aligned}\]
Then, we have
\[\begin{aligned}
     |H_{p,i}(y)|&\le C\ve_{p,i}|y|\quad\text{and}\quad \sup_{|y|\le\rho/\ve_{p,i}}\frac{p^m|H_{p,i}(y)|}{1+|y|^\tau} \le Cp^m\ve_{p,i}^{\tau}\to0,~~m=1,2,3.
\end{aligned}\]
Thus the opposite peak changes only the harmonic data by a term vanishing in the three weighted remainder norms. The decay in (P3) allows the localized arguments of the cited lemmas. The zero value and gradient at the origin exclude the dilation and translation kernel terms, and the annular estimates absorb the above harmonic error. Consequently, \cite[(5.19), (5.27) and (5.39)]{Ianni-Luo-Yan-2025} give
\[
 \sup_{|y|\le\rho/\ve_{p,i}}
 \frac{|v_{p,i}(y)|+|k_{p,i}(y)|+|q_{p,i}(y)|}
 {1+|y|^\tau}\le C,
\]
which proves \eqref{est-wpi-4}.

The coefficient equations follow by expanding $p\log(1+w_{p,i}/p)$ in those successive remainder equations. The coefficient at order $p^{-2}$ on the nonlinear side is
\[
e^U\left(w_1-Uw_0+\frac13U^3
+\frac12(w_0-U^2/2)^2\right),
\]
which gives \eqref{eq-w1-correct}. Local elliptic estimates for the remainder equation recover the stated $C^2$ bound, as in \cite[Proposition 5.1]{Ianni-Luo-Yan-2025}.
\end{proof}

\begin{lemma}\label{lem-mass-energy}
For every low-energy nodal family $u_p$, we have
\be\lab{energy}
p\int_{B_1}|u_p|^p\,\rd x
=16\pi\sqrt e\left(1-\frac{\log p}{p}\right)+O(p^{-1}),
\ee
and
\be\label{eq-energy-window-expansion}
p\int_{B_1}|\nabla u_p|^2\,\rd x
=16\pi e\left(1-\frac{2\log p}{p}\right)+O(p^{-1}).
\ee
In particular, $u_p\in\ER_p$ for all sufficiently large $p$.
\end{lemma}
\begin{proof}
A sharper expansion of $J_p(u_p)$ appears in \cite[Corollary~1.5]{Chen-Cheng-Zhao-2024}. We derive the estimates above directly to specify the energy window. Let $M_{p,i}=|u_p(x_{p,i})|$ and fix $\rho\in(0,\rho_0)$.
By (P5), the integrals outside the two fixed peak balls are $O(C^p/p^{p-1})$ after multiplication by $p$.
The change of variables $x=x_{p,i}+\ve_{p,i}y$ gives
\[
p\int_{B_\rho(x_{p,i})}|u_p|^p\,\rd x
=M_{p,i}\int_{B_{\rho/\ve_{p,i}}}
\left|1+\frac{w_{p,i}}p\right|^p\,\rd y.
\]
Properties (P4) and (P6) give \eqref{energy} after summation.
Testing the equation with $u_p$ and changing variables as above, we obtain
\[
p\int_{B_1}|\nabla u_p|^2\,\rd x
=\sum_{i=1}^2M_{p,i}^2
\int_{B_{\rho/\ve_{p,i}}}
\left|1+\frac{w_{p,i}}p\right|^{p+1}\,\rd y
+O(C^p/p^{p-1}).
\]
Now $M_{p,i}^2=e(1-2\log p/p)+O(p^{-1})$, so (P6) yields
\eqref{eq-energy-window-expansion}. In particular,
\[
16\pi e-p\int_{B_1}|u_p|^{p+1}\,\rd x
=32\pi e\frac{\log p}{p}+O(p^{-1}).
\]
This expansion implies both strict inequalities in \eqref{define-e_p} for sufficiently large $p$.
\end{proof}

\subsection{Limit eigenvalue problem}
\begin{lemma}\lab{lemma-3-10-1}
Let $\la\in[0,1]$, let $U$ be defined by \eqref{U}, and let $v\in C^2(\R^2)\setminus\{0\}$ solve
\[\begin{cases}
    -\Delta v=\la e^Uv \quad\text{in}~\R^2,\\
    \int_{\R^2}|\nabla v|^2 \rd x<\iy.
\end{cases}\]
Then either $\la=0$ or $\la=1$.
\end{lemma}
\begin{proof}
See \cite[p.~364]{DeMarchis-Grossi-Ianni-Pacella-2019} and \cite[Lemma~2.4]{Ianni-Luo-Yan-2025}.
\end{proof}

\begin{lemma}\label{lemma-limit-kernel}
Let $U$ be the function defined in \eqref{U}, and let $v\in C^2(\R^2)\setminus\{0\}$ solve
\[\begin{cases}
    -\Delta v= e^Uv \quad\text{in}~\R^2,\\
    \int_{\R^2}|\nabla v|^2 \rd x<\iy.
\end{cases}\]
Then
\[v\in\operatorname{span}\lbr{\f{x_1}{8+|x|^2},\f{x_2}{8+|x|^2},\f{8-|x|^2}{8+|x|^2}}.\]
\end{lemma}

\begin{proof}
See \cite[p.~364]{DeMarchis-Grossi-Ianni-Pacella-2019} and \cite[Lemma~2.5]{Ianni-Luo-Yan-2025}.
\end{proof}

\vskip0.2in

\section{Axial symmetry}\lab{sec-3}
We prove axial symmetry by rotating lines through the origin. The comparison argument follows the rotating-plane method used in \cite[Section 3]{Wei-Winter-2005}. Here the smallness of the difference quotient away from the two concentration points follows from \eqref{out-u}.

We first recall a stability property of a strictly concave radial function, see \cite[Lemma 4.2]{Ni-Takagi-1991}.
\begin{lemma} \lab{lemma-3-6-1}
Let $\phi\in C^2(\ov{B_a(0)})$ be radial, with $\phi(0)=0$ and $\phi''(r)<0$ for $0\le r\le a$. There exists $\delta>0$ such that every $\psi\in C^2(\ov{B_a(0)})$ satisfying
\[
    \nabla\psi(0)=0,\qquad
    \nm{\psi-\phi}_{C^2(\ov{B_a(0)})}<\delta
\]
has no critical point in $\ov{B_a(0)}\setminus\{0\}$.
\end{lemma}

\begin{lemma}\lab{lemma-3-6-2}
There exists $p_0>1$ such that, for every $p>p_0$ and every $u_p\in\ER_p\cap\CR_p$, the global maximum point and the global minimum point of $u_p$ are unique.
\end{lemma}
\begin{proof}
Suppose, for contradiction, that there exist $p_n\to+\infty$ and solutions $u_{p_n}\in\ER_{p_n}\cap\CR_{p_n}$, each with two distinct global maximum points $x_{p_n,1}$ and $\ti x_{p_n,1}$. This sequence satisfies the low-energy limit. We may therefore apply (P1)--(P6) after passing to a subsequence, and we write $p$ for $p_n$.

By (P3), (P5), and $u_p(x_{p,1})\to\sqrt e$, every global maximum point belongs to $B_\rho(x_{p,1})$ for large $p$, where $\rho\in(0,\rho_0)$ is fixed. Define
\[
    \eta_p=\frac{\ti x_{p,1}-x_{p,1}}{\ve_{p,1}}.
\]
Since $w_{p,1}(\eta_p)=0$, estimate \eqref{est-wpi-1}, with a fixed $\delta\in(0,1)$, gives
\[
    (4-\delta)\log|\eta_p|\le C_\delta.
\]
Thus $\eta_p$ is bounded. The convergence $w_{p,1}\to U$ in $C^2_{\loc}(\R^2)$ implies that every limit $\eta$ satisfies $U(\eta)=0$. Since $U$ attains its maximum only at the origin, $\eta_p\to0$.

The function $U$ satisfies $U(0)=0$ and $U''(r)<0$ on $[0,1]$. Applying Lemma \ref{lemma-3-6-1} to $w_{p,1}$ on $\ov{B_1(0)}$, we conclude that its only critical point there is the origin. Both $0$ and $\eta_p$ are critical points, so $\eta_p=0$ for large $p$, contrary to the choice of the two maximum points. Applying the same argument to $-u_p$, which belongs to $\ER_p\cap\CR_p$ and solves \eqref{equ-1}, proves uniqueness of the minimum point.
\end{proof}

\begin{proposition}\lab{lemma-1}
There exists $p_0>1$ such that, for every $p>p_0$ and every $u_p\in\ER_p\cap\CR_p$, the unique global maximum point $x_{p,1}$, the unique global minimum point $x_{p,2}$, and the origin are collinear. Moreover, $u_p$ is symmetric with respect to the line containing these points.
\end{proposition}
\begin{proof}
Consider a sequence of solutions $u_p\in\ER_p\cap\CR_p$ with $p\to+\infty$. As in the proof of Lemma~\ref{lemma-3-6-2}, this is a low-energy nodal family. Let $a=\sqrt{\sqrt5-2}$. By (P1) and the classification of critical points in Section~\ref{sec-2}, we may rotate and reflect coordinates so that
\begin{equation}\label{eq-assume}
    x_{p,1}=(t_p,s_{p,1}),\qquad
    x_{p,2}=(t_p,-s_{p,2}),\qquad
    t_p\ge0,\quad s_{p,1},s_{p,2}>0,
\end{equation}
where $t_p\to0$ and $s_{p,1},s_{p,2}\to a$. Indeed, we choose the line through the two extrema to be vertical and its signed distance from the origin to be nonnegative. The foot of the perpendicular from the origin lies between the extrema for large $p$, because their limits are antipodal.

Let $\theta_p\in(0,\pi/2]$ be the polar angle of $x_{p,1}$. For $0\le\theta\le\theta_p$, define
\[
    \Pi_\theta=\{x\in B_1(0):x_1\sin\theta-x_2\cos\theta=0\},\]
and
\[
    \Sigma_\theta=\{x\in B_1(0):x_1\sin\theta-x_2\cos\theta<0\}.
\]
We write $x^\theta$ for the reflection of $x$ across $\Pi_\theta$ and set
\[
    w_p^\theta(x)=u_p(x)-u_p(x^\theta),\qquad x\in\Sigma_\theta.
\]
The half-disks $\Sigma_\theta$ vary continuously, and $x_{p,1}\in\Sigma_\theta$ for $\theta<\theta_p$. The reflection difference satisfies
\begin{equation}\label{eq-reflection-equation}
\begin{cases}
    -\Delta w_p^\theta=c_p^\theta w_p^\theta&\text{in }\Sigma_\theta,\\
    w_p^\theta=0&\text{on }\pa\Sigma_\theta,
\end{cases}
\end{equation}
where the continuous difference quotient is defined at every point by
\begin{equation}\label{3-6-3}
    c_p^\theta(x)=p\int_0^1
    \left|\tau u_p(x)+(1-\tau)u_p(x^\theta)\right|^{p-1}\rd\tau.
\end{equation}
In particular,
\[
    0\le c_p^\theta(x)\le p\max\{|u_p(x)|,|u_p(x^\theta)|\}^{p-1}
    \le M_p:=p\nm{u_p}_{L^\infty(B_1(0))}^{p-1}.
\]

\smallskip
\noindent\textit{Step 1. Comparison at $\theta=0$.}
Choose $\rho\in(0,\rho_0)$ sufficiently small that $4\rho<a/2$. Since $s_{p,1}-s_{p,2}\to0$, reflection across $\Pi_0$ maps $B_{2\rho}(x_{p,1})$ into $B_{3\rho}(x_{p,2})$ for large $p$. By \eqref{3-6-2},
\[
    w_p^0\ge\frac2p\qquad\text{in }B_{2\rho}(x_{p,1}).
\]
Consequently, the set $\{x\in\Sigma_0:w_p^0(x)<0\}$ is disjoint from this ball. Every point of that set, together with its reflection, stays at least a fixed positive distance from both concentration points. More precisely, $x_2>0$ separates $x$ from $x_{p,2}$ and $x^0$ from $x_{p,1}$, while
\[
    |x^0-x_{p,2}|\ge |x-x_{p,1}|-|s_{p,1}-s_{p,2}|\ge\rho
\]
for large $p$. By \eqref{out-u}, both $|u_p(x)|$ and $|u_p(x^0)|$ are bounded by $C/p$ on this set, with $C$ independent of $p$. Hence
\[
    c_p^0(x)\le p(C/p)^{p-1}=o(1)
    \qquad\text{on }\{w_p^0<0\}.
\]
Let $v_p=\max\{-w_p^0,0\}\in H_0^1(\Sigma_0)$. Testing \eqref{eq-reflection-equation} with $-v_p$ and using the Poincar\'e inequality, we obtain
\[
    \int_{\Sigma_0}|\nabla v_p|^2\rd x
    =\int_{\Sigma_0}c_p^0v_p^2\rd x
    \le\frac{o(1)}{\lambda_1(\Sigma_0)}
        \int_{\Sigma_0}|\nabla v_p|^2\rd x,
\]
where $\lambda_1(\Sigma_0)>0$ is the first Dirichlet eigenvalue of $-\Delta$. Thus $v_p=0$ for large $p$. Since $w_p^0\ge2/p$ near $x_{p,1}$, the strong maximum principle gives
\begin{equation}\label{3-6-6}
    w_p^0>0\qquad\text{in }\Sigma_0.
\end{equation}

\smallskip
\noindent\textit{Step 2. Continuation of the comparison.}
Fix a sufficiently large $p$. We first show that strict positivity at an angle $\bar\theta<\theta_p$ persists for nearby angles. Suppose that $w_p^{\bar\theta}>0$ in $\Sigma_{\bar\theta}$. Let $C_S$ be a constant in the Sobolev inequality
\[
    \nm{v}_{L^4(B_1(0))}^2\le C_S\int_{B_1(0)}|\nabla v|^2\rd x,
    \qquad v\in H_0^1(B_1(0)).
\]
Choose a compact set $K\subset\Sigma_{\bar\theta}$ so that its complement has sufficiently small area. For $\theta$ close to $\bar\theta$, we have $K\subset\Sigma_\theta$, $w_p^\theta>0$ on $K$, and
\[
    C_SM_p|\Sigma_\theta\setminus K|^{1/2}<1.
\]
Here $p$ is fixed, so $M_p$ is fixed. The negative part $v=\max\{-w_p^\theta,0\}$ vanishes on $K$ and belongs to $H_0^1(\Sigma_\theta)$. Extend it by zero to $B_1(0)$. Testing \eqref{eq-reflection-equation} and applying H\"older's inequality, we obtain
\[
\begin{aligned}
    \int_{B_1(0)}|\nabla v|^2\rd x
    &\le M_p\int_{\Sigma_\theta\setminus K}v^2\rd x\\
    &\le C_SM_p|\Sigma_\theta\setminus K|^{1/2}
        \int_{B_1(0)}|\nabla v|^2\rd x.
\end{aligned}
\]
Therefore $w_p^\theta\ge0$. Since $x_{p,1}\in\Sigma_\theta$ for $\theta<\theta_p$ and the maximum point is unique,
\[
    w_p^\theta(x_{p,1})=u_p(x_{p,1})-u_p(x_{p,1}^\theta)>0.
\]
The strong maximum principle yields $w_p^\theta>0$ in $\Sigma_\theta$.

Define
\[
    \theta_*:=\sup\{\alpha\in[0,\theta_p):
    w_p^\theta>0\text{ in }\Sigma_\theta
    \text{ for every }\theta\in[0,\alpha]\}.
\]
By \eqref{3-6-6} and the continuation argument, $\theta_*>0$. If $\theta_*<\theta_p$, continuity gives $w_p^{\theta_*}\ge0$, and uniqueness of the maximum point implies $w_p^{\theta_*}(x_{p,1})>0$. The strong maximum principle and the same continuation argument extend the comparison beyond $\theta_*$, a contradiction. Hence
\begin{equation}\label{eq-rotating-comparison}
    w_p^\theta>0\text{ in }\Sigma_\theta
    \quad(0\le\theta<\theta_p),\qquad
    w_p^{\theta_p}\ge0\text{ in }\Sigma_{\theta_p}.
\end{equation}

\smallskip
\noindent\textit{Step 3. The limiting line.}
If $w_p^{\theta_p}$ is not identically zero, it is positive in $\Sigma_{\theta_p}$. The point $x_{p,1}$ lies on the straight part of $\pa\Sigma_{\theta_p}$ and strictly inside $B_1(0)$. The Hopf boundary lemma therefore gives a strictly positive inward normal derivative of $w_p^{\theta_p}$ at $x_{p,1}$. On the other hand, $x_{p,1}$ is a critical point of $u_p$ and is fixed by the reflection, so
\[
    \nabla w_p^{\theta_p}(x_{p,1})=0.
\]
This contradiction proves $w_p^{\theta_p}\equiv0$. Thus $u_p$ is symmetric across $\Pi_{\theta_p}$. Its unique minimum point must also be fixed by this reflection. Hence $\Pi_{\theta_p}$ contains both extrema and the origin. Since the line through the two points in \eqref{eq-assume} is $\{x_1=t_p\}$, we conclude that $t_p=0$ and $\theta_p=\pi/2$.

The argument includes $t_p=0$ from the outset and proves the stated properties for every sufficiently large member of the chosen sequence. If no uniform $p_0$ existed, we could choose a sequence violating the proposition and obtain a contradiction by this argument.
\end{proof}

\begin{corollary}\label{cor-angular-monotonicity}
For sufficiently large $p$, rotate a solution $u_p\in\ER_p\cap\CR_p$ so that its positive maximum point lies on the positive $x_1$-axis. In polar coordinates,
\begin{equation}\label{eq-angular-monotonicity}
    \pa_\theta u_p(r,\theta)<0
    \qquad\text{for }0<r<1,\quad0<\theta<\pi.
\end{equation}
\end{corollary}
\begin{proof}
First use the coordinates in the proof of Proposition \ref{lemma-1}, in which the positive maximum lies on the positive $x_2$-axis. Write $u_p=u_p(r,\varphi)$. By \eqref{eq-rotating-comparison}, $w_p^\theta>0$ for $0\le\theta<\pi/2$. For $0<r<1$, apply the Hopf boundary lemma at the points $r(\cos\theta,\sin\theta)$ and $-r(\cos\theta,\sin\theta)$ on the straight boundary of $\Sigma_\theta$. The inward unit normal is $(-\sin\theta,\cos\theta)$, so the reflection formula gives
\[
    \pa_\varphi u_p(r,\theta)>0,\qquad
    \pa_\varphi u_p(r,\theta+\pi)<0
    \qquad(0\le\theta<\pi/2).
\]
Symmetry across the $x_2$-axis gives
\[
    \pa_\varphi u_p(r,\pi-\varphi)
    =-\pa_\varphi u_p(r,\varphi).
\]
Combining these relations yields $\pa_\varphi u_p<0$ for $\pi/2<\varphi<3\pi/2$. Finally, rotate the coordinates by $\pi/2$ so that the positive maximum lies on the positive $x_1$-axis. In the new coordinates, $\varphi=\theta+\pi/2$, which proves \eqref{eq-angular-monotonicity}.
\end{proof}

\vskip0.2in
\section{Morse index in the symmetric space}\lab{sec-4}
By Proposition \ref{lemma-1}, every low-energy nodal solution is symmetric with respect to a line through the origin for sufficiently large $p$. We rotate this line onto the $x_1$-axis and work in the corresponding symmetric space.
\begin{definition}
A function $u\in H_0^1(B_1(0))$ is symmetric with respect to the $x_1$-axis if
    \[u(x_1,x_2)=u(x_1,-x_2).\]
We define 
\be\lab{HRe} \HR_e:=\{u\in H_0^1(B_1(0)):~u(x_1,x_2)=u(x_1,-x_2) \}. \ee 
The space $\HR_e$ is a closed subspace of $H_0^1(B_1(0))$, with inner product
    \[\abr{u,v}=\int_{B_1(0)}\nabla u\cdot\nabla v \rd x.\]
\end{definition}

We denote by $m_e(u_p)$ the Morse index of the symmetric solution $u_p\in\CR_p$ in $\HR_e$, defined as the maximal dimension of a subspace on which the quadratic form $J_p''(u_p)$ is negative definite.
The variational characterization gives 
\be  m_e(u_p)=\sharp\{j\in\N:\la_{p,j}<1 \} \ee
where $\la_{p,1}<\la_{p,2}\le \la_{p,3}\le \cdots$ denotes the sequence of eigenvalues of the problem
\be \begin{cases}
    -\Delta v=\la p|u_p|^{p-1}v \quad\text{in}~B_1(0),\\
    v=0 \quad\text{on}~\pa B_1(0),\\
    v\in\HR_e.
\end{cases}\ee
We prove that every even low-energy nodal solution has Morse index two and a trivial linearized kernel in $\HR_e$.

\begin{proposition}\lab{thm-index}
There exists $p_0>1$ such that, for every $p>p_0$ and every $u_p\in\ER_p\cap\CR_p\cap\HR_e$, 
    \[m_e(u_p)=2,\]
and $u_p$ is nondegenerate in $\HR_e$.
\end{proposition}

Fix a sequence of solutions $u_p\in\ER_p\cap\CR_p\cap\HR_e$ with $p\to+\iy$. It suffices to prove the asserted index and nondegeneracy along every such sequence.
The blow-up points $x_{p,1}$ and $x_{p,2}$ both lie on the $x_1$-axis.
Replacing $u_p(x)$ by $u_p(-x)$ if necessary, we assume
\be\lab{assumption-1}  (x_{\iy,1},x_{\iy,2})=(X_1,X_2),\ee 
where $X_1,X_2$ are given by \eqref{X}.
The map $v(x)\mapsto v(-x)$ preserves $\HR_e$ and conjugates the two linearized problems, so the opposite orientation has the same spectrum.
We use the spectral approach developed in \cite{DeMarchis-Grossi-Ianni-Pacella-2019,Ianni-Luo-Yan-2025}. Set $\sigma_j=(-1)^{j-1}$. By (P3), $\sigma_j u_p$ is positive and satisfies the positive Lane--Emden equation in a fixed ball about $x_{p,j}$. The rescaled potential therefore has the same form at both peaks. The signs $\sigma_j$ enter the outer expansion of $u_p$ and the Pohozaev identities; they will be retained in the computation below.

Let $v_{p,l}$ be an eigenfunction associated with $\la_{p,l}$, so that
\be\lab{equ-v}\begin{cases}
    -\Delta v_{p,l}=\la_{p,l} p|u_p|^{p-1}v_{p,l} \quad\text{in}~B_1(0),\\
    v_{p,l}=0 \quad\text{on}~\pa B_1(0),\\
    v_{p,l}\in\HR_e.
\end{cases}\ee
Choose an orthogonal basis in each eigenspace. Eigenfunctions belonging to distinct eigenvalues are already orthogonal, and hence
\be \int_{B_1(0)}\nabla v_{p,l}\cdot\nabla v_{p,l'}\rd x=0  \quad\text{if}~l\neq l'. \ee
Normalize each eigenfunction by
\be\lab{v-1} \nm{v_{p,l}}_{L^\iy(B_1(0))}=1. \ee 

\vskip0.1in
\subsection{Estimates for the eigenfunctions}\ \\ \indent
For $l\in\N$ and $j=1,2$, define
\be\lab{ti-v} \ti v_{p,l,j}(x)=v_{p,l}(x_{p,j}+\ve_{p,j}x),\quad x\in\Omega_{p,j}. \ee 
By \eqref{equ-v} and \eqref{v-1}, it follows that
\be\lab{equ-3-10-1}\begin{cases}
    -\Delta \ti v_{p,l,j}=\la_{p,l}\wt V_{p,j}\ti v_{p,l,j} \quad\text{in}~\Omega_{p,j},\\
    \ti v_{p,l,j}=0 \quad\text{on}~\pa \Omega_{p,j},\\
    \ti v_{p,l,j}(x)=\ti v_{p,l,j}(\bar x),\\
    \nm{\ti v_{p,l,j}}_{L^\iy(\Omega_{p,j})}=1,
\end{cases}\ee
where $\wt V_{p,j}=\abs{1+\f{w_{p,j}}{p}}^{p-1}$ with $w_{p,j}$ in \eqref{w_pi}, and we denote
    \[\bar x=(x_1,-x_2)\quad\text{for}~x=(x_1,x_2)\in\R^2.\]
Whenever the eigenvalues under consideration are bounded, we pass to a subsequence and write
\[\lim_{p\to\iy}\la_{p,l}=\la_{\iy,l}.\]
The values relevant to the Morse index are $\la_{\iy,l}=0$ and $\la_{\iy,l}=1$. 

\begin{lemma}\lab{lemma-2}
Let $l\in\N$ and suppose that $\la_{p,l}\to\la_{\iy,l}$. After passing to a subsequence, the following conclusions hold.
\begin{itemize}
\item[(1)] If $\la_{\iy,l}=0$, then for any $j=1,2$, there exists $c_{l,j}\in\R$ such that
    \[\ti v_{p,l,j}\to c_{l,j} \quad\text{in}~C_\loc^1(\R^2),~\text{as}~p\to\iy,\]
with $(c_{l,1},c_{l,2})\neq (0,0)$. Moreover, if $\la_{\iy,l}=\la_{\iy,l'}=0$ for $l\neq l'$, then we have
\be  \sum_{j=1}^2c_{l,j}c_{l',j}=0. \ee
\item[(2)] If $\la_{\iy,l}=1$, then for any $j=1,2$, there exists $(a_{l,j},b_{l,j})\in\R^2$ such that 
\be\lab{in-vplj} \ti v_{p,l,j}\to V_{l,j}=a_{l,j}\f{x_1}{8+|x|^2}+b_{l,j}\f{8-|x|^2}{8+|x|^2} \quad\text{in}~C_\loc^1(\R^2),~\text{as}~p\to\iy,\ee 
and either $(a_{l,1},b_{l,1})\neq (0,0)$ or  $(a_{l,2},b_{l,2})\neq (0,0)$. Moreover, if $\la_{\iy,l}=\la_{\iy,l'}=1$ for $l\neq l'$, then we have
\be \sum_{j=1}^2\bigl(a_{l,j}a_{l',j}+32b_{l,j}b_{l',j}\bigr)=0.\ee
\end{itemize}
\end{lemma}
\begin{proof}
Local elliptic estimates give, after extraction, bounded limits $V_{l,j}$ satisfying
\[-\Delta V_{l,j}=\la_{\iy,l}e^UV_{l,j}\quad\text{in }\R^2.\]
We first show that the two local limits $V_{l,j}$ cannot both vanish. This argument applies whenever $\la_{p,l}$ is bounded. Suppose that both limits are zero. Green's representation formula and the concentration estimates in (P3)--(P5) imply
\[v_{p,l}\to0\quad\text{uniformly on }\ov{B_1(0)}\setminus\bigcup_{j=1}^2B_d(x_{p,j})\]
for every sufficiently small fixed $d>0$. Fix $\delta\in(0,1)$. In the rescaled annulus $R<|y|<d/\ve_{p,j}$, (P3) gives
\[|\Delta\ti v_{p,l,j}(y)|\le C|y|^{-4+\delta}.\]
Comparison with the radial solution of $-\Delta h=C|y|^{-4+\delta}$ on this annulus, with boundary values equal to the maximum of $|\ti v_{p,l,j}|$ on the two boundary circles, yields
\[\limsup_{p\to\iy}\sup_{R<|y|<d/\ve_{p,j}}|\ti v_{p,l,j}(y)|\le CR^{-2+\delta}.\]
The boundary values tend to zero for fixed $R$, by local convergence at $|y|=R$ and the preceding outer estimate at $|y|=d/\ve_{p,j}$. Letting $R\to\iy$ contradicts \eqref{v-1}. Thus at least one local limit is nonzero.

Testing the limit equation with logarithmic cutoffs shows that each bounded limit has finite Dirichlet energy. If $\la_{\iy,l}=0$, Liouville's theorem makes each $V_{l,j}$ constant. If $\la_{\iy,l}=1$, the classification of bounded solutions in Lemma \ref{lemma-limit-kernel} gives the translation and dilation terms in \eqref{in-vplj}. Evenness in $x_2$ removes the second translation term.

For $l\ne l'$, the eigenvalue equation and Dirichlet orthogonality give
\[p\int_{B_1(0)}|u_p|^{p-1}v_{p,l}v_{p,l'}\rd x=0.\]
The integral outside the two fixed peak neighborhoods tends to zero. Rescaling the two remaining integrals and using (P3), we obtain
\[\sum_{j=1}^2\int_{\R^2}e^UV_{l,j}V_{l',j}\rd y=0.\]
For constant limits this gives the assertion in (1). For the limits in (2), the mixed translation--dilation integrals vanish by oddness, while
\[\int_{\R^2}e^U\frac{y_1^2}{(8+|y|^2)^2}\rd y=\frac\pi{12},\qquad
      \int_{\R^2}e^U\left(\frac{8-|y|^2}{8+|y|^2}\right)^2\rd y=\frac{8\pi}{3}.\]
This proves the stated global orthogonality relation.
\end{proof}

We next expand the eigenfunctions away from the concentration points.

\begin{lemma}\lab{lemma-3}
For fixed $l\in\N$ with bounded $\la_{p,l}$ and every sufficiently small fixed $d>0$,
\be\lab{out-vpl-1} v_{p,l}(x)=\la_{p,l}p \sum_{j=1}^2 A_{p,l,j}G(x_{p,j},x)+\la_{p,l}p\sum_{j=1}^2 B_{p,l,j}\pa_{x_1}G(x_{p,j},x)+o(\la_{p,l}\ve_p), \ee
in $C^1(\ov{B_1(0)}\setminus(B_{2d}(X_1)\cup B_{2d}(X_2)))$, where $\ve_p$ is defined by \eqref{ve_p} and 
\be\lab{A_plj} A_{p,l,j}=\int_{B_d(x_{p,j})}|u_p|^{p-1}v_{p,l}\rd x, \ee 
\be\lab{B_plj} B_{p,l,j}=\int_{B_d(x_{p,j})}(x-x_{p,j})_1|u_p|^{p-1}v_{p,l}\rd x. \ee 
\end{lemma}
\begin{proof} 
Apply Green's representation formula to $p|u_p|^{p-1}v_{p,l}$ and expand the Green function to first order in each $B_d(x_{p,j})$. The moment in the $x_2$ direction vanishes by symmetry. The second-order remainder is $o(\ve_p)$ after multiplication by $p$, by \eqref{est-wpi-3}; the contribution outside the two balls is smaller than any fixed power of $\ve_p$. This is the expansion in \cite[Lemma 4.2]{Chen-Cheng-Zhao-2024}, applied to the present eigenfunctions. Its proof uses the nonnegative potential $|u_p|^{p-1}$ and thus applies to both peak signs.
\end{proof}

For $j=1,2$, define the bilinear forms
\be\lab{Pj} P_j(u,v)=-2r\int_{\pa B_r(x_{p,j})} \pa_\nu u\pa_\nu v \rd\sigma +r\int_{\pa B_r(x_{p,j})}\nabla u\cdot\nabla v \rd\sigma, \ee
and
\be\lab{Qj} Q_j(u,v)=-\int_{\pa B_r(x_{p,j})} \pa_\nu u\pa_{x_1} v+\pa_\nu v\pa_{x_1} u \rd\sigma +\int_{\pa B_r(x_{p,j})}(\nabla u\cdot\nabla v)\nu_1 \rd\sigma, \ee
where $u$ and $v$ are $C^2$ in a neighborhood of $\partial B_r(x_{p,j})$, and $r>0$ is small enough that the balls $B_{2r}(x_{p,1})$ and $B_{2r}(x_{p,2})$ are disjoint and contained in $B_1(0)$.
For the Green function and its derivatives, these forms satisfy the identities below; see \cite{Cao-Guo-Peng-Yan-2019,Grossi-Ianni-Luo-Yan-2022}.
\begin{lemma}\lab{lemma-3-11-1}
We have
\[P_j\Big(G(x_{p,s},x),G(x_{p,m},x)\Big)=\begin{cases} 
    -\f{1}{2\pi} \quad &\text{for}~s=m=j,\\
    0 \quad &\text{for}~s\neq j~\text{or}~m\neq j,
\end{cases}\]
\[P_j\Big(G(x_{p,s},x),\pa_1G(x_{p,m},x)\Big)=\begin{cases} 
    \f{1}{2}\pa_1R(x_{p,j}) \quad &\text{for}~s=m=j,\\
    D_1G(x_{p,s},x_{p,j}) \quad &\text{for}~m=j,s\neq j,\\
    0 \quad &\text{for}~m\neq j,
\end{cases}\]
\[Q_j\Big(G(x_{p,s},x),G(x_{p,m},x)\Big)=\begin{cases} 
    \pa_{x_1}R(x_{p,j}) \quad &\text{for}~s=m=j,\\
    D_1G(x_{p,m},x_{p,j})\quad &\text{for}~s=j,m\neq j,\\
    D_1G(x_{p,s},x_{p,j})\quad &\text{for}~s\neq j,m=j,\\
    0\quad &\text{for}~s\neq j,m\neq j,\\
\end{cases}\]
and
\[Q_j\Big(G(x_{p,s},x),\pa_1G(x_{p,m},x)\Big)=\begin{cases} 
    \f12\pa_{x_1x_1}^2R(x_{p,j}) \quad &\text{for}~s=m=j,\\
    D_1\pa_1G(x_{p,m},x_{p,j})\quad &\text{for}~s=j,m\neq j,\\
    D_{11}^2G(x_{p,s},x_{p,j})\quad &\text{for}~s\neq j,m=j,\\
    0\quad &\text{for}~s\neq j,m\neq j,\\
\end{cases}\]
where $\pa_1G(x,y)=\pa_{x_1}G(x,y)$ and $D_1G(x,y)=\pa_{y_1}G(x,y)$.
\end{lemma}

We have also the following local Pohozaev identities.
\begin{lemma}
For each $j=1,2$,
\be\lab{Pj-1} \begin{aligned}
    P_j(u_p,v_{p,l})&=\int_{\pa B_r(x_{p,j})}|u_p|^{p-1}u_pv_{p,l} (x-x_{p,j})\cdot\nu\rd\sigma -2\int_{B_r(x_{p,j})}|u_p|^{p-1}u_pv_{p,l}\rd x\\
    &\quad +(\la_{p,l}-1)p\int_{B_r(x_{p,j})}|u_p|^{p-1}v_{p,l}(x-x_{p,j})\cdot\nabla u_p\rd x, 
\end{aligned}\ee
and
\be\lab{Qj-1} \begin{aligned}
    Q_j(u_p,v_{p,l})&=\int_{\pa B_r(x_{p,j})}|u_p|^{p-1}u_pv_{p,l} \nu_1\rd\sigma +(\la_{p,l}-1)p\int_{B_r(x_{p,j})}|u_p|^{p-1}v_{p,l}\pa_{x_1}u_p\rd x.
\end{aligned}\ee
\end{lemma}
\begin{proof}
Set $f(t)=|t|^{p-1}t$, so that $f'(t)=p|t|^{p-1}$, and write $z=x-x_{p,j}$. The divergence of
\[(\nabla u_p\cdot\nabla v_{p,l})z-(z\cdot\nabla u_p)\nabla v_{p,l}-(z\cdot\nabla v_{p,l})\nabla u_p\]
is
\[\operatorname{div}(zf(u_p)v_{p,l})-2f(u_p)v_{p,l}
      +(\la_{p,l}-1)f'(u_p)v_{p,l}z\cdot\nabla u_p.\]
On $\partial B_r(x_{p,j})$, the normal component of this vector field is the integrand defining $P_j$. The divergence theorem gives \eqref{Pj-1}. Replacing $z$ by the constant vector $(1,0)$ removes the term $-2f(u_p)v_{p,l}$ and gives \eqref{Qj-1}.
\end{proof}

For eigenvalues converging to $1$, the coefficient of the radial kernel has the following consequences.
\begin{lemma}\lab{lemma-4}
Let $l\in\N$. If $\la_{\iy,l}=1$, then we have
\be\lab{3-11-5} A_{p,l,j}=-\f{8\pi}{p^2}(b_{l,j}+o(1)), \ee 
and
\be\lab{out-vpl-2} v_{p,l}=-\f{8\pi}{p}\sum_{j=1}^2 b_{l,j}G(x_{p,j},x)+o\sbr{\f1p},\quad \text{in}~C_\loc^1(\ov{B_1(0)}\setminus\{X_1,X_2\}), \ee
where $b_{l,j}$ are constants in \eqref{in-vplj}. Moreover, if $(b_{l,1},b_{l,2})\neq(0,0)$, then it holds
\be \la_{p,l}=1+\f{6}{p}+o\sbr{\f1p}. \ee 
\end{lemma}
\begin{proof}
Write $M_{p,j}=u_p(x_{p,j})$ and
\[F_{p,l,j}=\int_{B_d(x_{p,j})}|u_p|^{p-1}u_pv_{p,l}\rd x.\]
Since $\int_{\R^2}e^UV_{l,j}=0$, rescaling \eqref{A_plj} first gives $pA_{p,l,j}=o(1)$. Also $B_{p,l,j}=O(\ve_p/p)$ by (P3). Lemma \ref{lemma-3} therefore gives $v_{p,l}=o(1)$ in $C^1$ away from the two peaks. Green's identity on $B_d(x_{p,j})$ yields
\[(p\la_{p,l}-1)F_{p,l,j}
      =\int_{\partial B_d(x_{p,j})}\bigl(v_{p,l}\partial_\nu u_p-u_p\partial_\nu v_{p,l}\bigr)\rd\sigma,
      \qquad F_{p,l,j}=o(p^{-2}).\]
On the other hand, the rescaling gives the exact relation
\[\frac{F_{p,l,j}}{M_{p,j}}=A_{p,l,j}
      +\frac1{p^2}\int_{B_{d/\ve_{p,j}}(0)}w_{p,j}
       \left|1+\frac{w_{p,j}}p\right|^{p-1}\ti v_{p,l,j}\rd y.\]
By (P3) and Lemma \ref{lemma-2}, the last integral converges to
\[\int_{\R^2}Ue^UV_{l,j}\rd y=8\pi b_{l,j}.\]
This proves \eqref{3-11-5}, and Lemma \ref{lemma-3} gives \eqref{out-vpl-2}.

Suppose that $b_{l,j}\ne0$ for some $j$. The improved outer estimate $v_{p,l}=O(p^{-1})$ gives $F_{p,l,j}=O(p^{-3})$. Expanding $P_j(u_p,v_{p,l})$ with Lemma \ref{lemma-3-11-1} gives
\[P_j(u_p,v_{p,l})
      =-\frac{\la_{p,l}pC_{p,j}}{2\pi}A_{p,l,j}+o(p^{-2})
      =\frac{32\pi\sigma_j\sqrt e}{p^2}b_{l,j}+o(p^{-2}).\]
The boundary integral in \eqref{Pj-1} is $o(p^{-2})$. Rescaling the last integral there gives
\[\begin{aligned}
      &p\int_{B_d(x_{p,j})}|u_p|^{p-1}v_{p,l}(x-x_{p,j})\cdot\nabla u_p\rd x\\
      &\quad=\frac{M_{p,j}}p\left(\int_{\R^2}e^UV_{l,j}y\cdot\nabla U\rd y+o(1)\right)\\
      &\quad=\frac{\sigma_j\sqrt e}{p}\left(\frac{16\pi}{3}b_{l,j}+o(1)\right).
    \end{aligned}\]
For this passage to the limit, (P3) and interior gradient estimates give
$|\nabla w_{p,j}(y)|\le C\log(2+|y|)/(1+|y|)$ on the smaller rescaled ball, which provides an integrable bound. Substitution into \eqref{Pj-1} gives
\[\frac{32\pi\sigma_j\sqrt e}{p^2}b_{l,j}+o(p^{-2})
      =(\la_{p,l}-1)\frac{\sigma_j\sqrt e}{p}
        \left(\frac{16\pi}{3}b_{l,j}+o(1)\right).\]
Division by the nonzero limiting coefficient proves $\la_{p,l}=1+6/p+o(p^{-1})$. This calculation is local and the factor $\sigma_j$ cancels, as in the positive case treated in \cite[Proposition 3.5]{Ianni-Luo-Yan-2025}.
\end{proof}

\vskip0.1in
\subsection{Computation of the Morse index}\ \\ \indent
Let $r\in(0,\rho_0)$ be such that
\[B_{4r}(X_1)\subset B_1(0),\quad B_{4r}(X_2)\subset B_1(0),\quad B_{3r}(X_1)\cap B_{3r}(X_2)=\emptyset.\]
For each $j=1,2$, we set $\phi_j(x)=\phi(x-X_j)$, where $\phi\in C_c^\iy(\R^2)$ is radial and satisfies $\phi=1$ in $B_{2r}(0)$, $\phi=0$ in $\R^2\setminus B_{3r}(0)$, $\phi\in(0,1)$ in $B_{3r}(0)\setminus B_{2r}(0)$. Since the centers lie on the $x_1$-axis, both cutoffs belong to $\HR_e$.
We define
\be u_{p,j}=\phi_j u_p,\quad \psi_{p,j}=\phi_j\pa_{x_1}u_p. \ee 
Then $u_{p,j}\in \HR_e$ and $\psi_{p,j}\in\HR_e$ for each $j=1,2$.
Testing the equation for $u_p$ with $\phi_j^2u_p$ and the differentiated equation with $\phi_j^2\pa_{x_1}u_p$ gives 
\be\lab{3-10-1} \int_{B_1(0)}|\nabla u_{p,j}|^2\rd x=\int_{B_1(0)}|\nabla \phi_j|^2u_p^2+\phi_j^2|u_p|^{p+1}\rd x, \ee
\be\lab{3-10-2} \int_{B_1(0)}|\nabla\psi_{p,j}|^2\rd x=\int_{B_1(0)}|\nabla\phi_j|^2|\pa_{x_1}u_p|^2+p\phi_j^2|u_p|^{p-1}|\pa_{x_1}u_p|^2\rd x, \ee 
and
\be\lab{3-10-3} \int_{B_1(0)}\nabla u_{p,j}\cdot\nabla\psi_{p,j}\rd x=\int_{B_1(0)}|\nabla\phi_j|^2u_p\pa_{x_1}u_p+\f{p+1}{2}\phi_j^2|u_p|^{p-1}u_p\pa_{x_1}u_p\rd x, \ee 
for each $j=1,2$.

We first use $u_{p,1},u_{p,2}$ as test functions to estimate $\la_{p,1}$ and $\la_{p,2}$.
\begin{lemma}
We have
\be\lab{3-10-4} 0<\la_{p,1}<\la_{p,2}\le\f1p+O\sbr{\f{1}{p^2}}. \ee 
\end{lemma}
\begin{proof}
Let $E_p=\mathrm{span}\{u_{p,1},u_{p,2}\}$. The supports are disjoint and both functions are nonzero, so $\dim(E_p)=2$. The min--max principle gives
    \[\la_{p,1}=\min_{f\in\HR_e\setminus\{0\}}\f{\int_{B_1(0)}|\nabla f|^2\rd x}{p\int_{B_1(0)}|u_p|^{p-1}f^2\rd x}>0,\]
and 
    \[\la_{p,2}=\min_{\substack{E\subset\HR_e\\ \dim(E)=2}}\max_{f\in E\setminus\{0\}}\f{\int_{B_1(0)}|\nabla f|^2\rd x}{p\int_{B_1(0)}|u_p|^{p-1}f^2\rd x}\le \max_{f\in E_p\setminus\{0\}}\f{\int_{B_1(0)}|\nabla f|^2\rd x}{p\int_{B_1(0)}|u_p|^{p-1}f^2\rd x}. \]
For $f=t_1u_{p,1}+t_2u_{p,2}$ with $(t_1,t_2)\neq (0,0)$, disjointness of the supports gives
\[\begin{aligned}
    \f{\int_{B_1(0)}|\nabla f|^2\rd x}{p\int_{B_1(0)}|u_p|^{p-1}f^2\rd x} &=\f{\sum_{j=1}^2t_j^2\int_{B_1(0)}|\nabla u_{p,j}|^2\rd x}{p\sum_{j=1}^2t_j^2\int_{B_1(0)}|u_p|^{p-1}u_{p,j}^2\rd x}\\
    &\le \max_{j=1,2}\f{\int_{B_1(0)}|\nabla u_{p,j}|^2\rd x}{p\int_{B_1(0)}\phi_j^2|u_p|^{p+1}\rd x}\\
    &=\f1p+\max_{j=1,2}\f{\int_{B_1(0)}|\nabla \phi_j|^2u_p^2\rd x}{p\int_{B_1(0)}\phi_j^2|u_p|^{p+1}\rd x},
\end{aligned}\]
where we have used \eqref{3-10-1} in the last equality.
On the other hand, we have
    \[p\int_{B_1(0)}\phi_j^2|u_p|^{p+1}\rd x\ge p\int_{B_r(X_j)}|u_p|^{p+1}\rd x\ge c_0,\]
and by \eqref{out-u},
    \[\int_{B_1(0)}|\nabla \phi_j|^2u_p^2\rd x\le C\int_{B_{3r}(X_j)\setminus B_{2r}(X_j)}u_p^2\rd x=O\sbr{\f{1}{p^2}},\]
so that we deduce \eqref{3-10-4}.
\end{proof}

Next we use $u_{p,1},u_{p,2},\psi_{p,1},\psi_{p,2}$ as test functions to estimate $\la_{p,3}$ and $\la_{p,4}$.
\begin{lemma}
The eigenvalues satisfy $\la_{\iy,3}=\la_{\iy,4}=1$ and
\be\lab{la_34} \la_{p,3}\le \la_{p,4}\le 1+O(\ve_p^2),\ee
where $\ve_p$ is defined by \eqref{ve_p}.
\end{lemma}
\begin{proof}
We first prove $\la_{\iy,4}\ge\la_{\iy,3}\ge1$.  Suppose, for contradiction, that $\la_{\iy,3}<1$. Then for each $j=1,2$, since $\wt V_{p,j}$ is bounded in $B_R(0)$ by \eqref{est-wpi-1}, local elliptic estimates yield convergence of the rescaled eigenfunctions to bounded solutions $V_{3,j}$ of
\[-\Delta V_{3,j}=\la_{\iy,3}e^UV_{3,j} \quad\text{in}~\R^2,\quad \nm{V_{3,j}}_{L^\iy(\R^2)}\le 1.\]
At least one $V_{3,j}$ is nonzero by the localization argument in Lemma \ref{lemma-2}. Lemma \ref{lemma-3-10-1} then gives $\la_{\iy,3}=0$. Lemma \ref{lemma-2} implies $V_{3,j}=c_{3,j}$ and $(c_{3,1},c_{3,2})\neq(0,0)$. 
By \eqref{3-10-4}, $\la_{\iy,1}=\la_{\iy,2}=0$, and Lemma \ref{lemma-2} implies that there exist $(c_{1,1},c_{1,2})\neq (0,0)$ and $(c_{2,1},c_{2,2})\neq (0,0)$ such that
    \[(c_{l,1},c_{l,2})\perp(c_{l',1},c_{l',2}),\quad\text{for any}~l\neq l'~\text{and}~l,l'\in\{1,2,3\}.\]
Three nonzero vectors in $\R^2$ cannot be pairwise orthogonal. This contradiction gives $\la_{\iy,4}\ge\la_{\iy,3}\ge1$. The upper bound \eqref{la_34}, proved next, therefore gives $\la_{\iy,3}=\la_{\iy,4}=1$.

Let
\[W_p=\mathrm{span}\{u_{p,1},u_{p,2},\psi_{p,1},\psi_{p,2}\}\subset\HR_e.\]
At each peak, $u_p$ is nonzero whereas $\pa_{x_1}u_p$ vanishes at the peak and is not identically zero in its neighborhood. Together with the disjoint supports, this gives $\dim(W_p)=4$. The min--max principle yields
    \[\la_{p,4}\le \max_{f\in W_p\setminus\{0\}}\f{\int_{B_1(0)}|\nabla f|^2\rd x}{p\int_{B_1(0)}|u_p|^{p-1}f^2\rd x}. \]
For any $f\in W_p\setminus\{0\}$, we can write
    \[f=f_1+f_2\quad\text{with}\quad f_1=\sum_{j=1}^2\al_{p,j}u_{p,j},~f_2=\sum_{j=1}^2\beta_{p,j}\psi_{p,j}\]
where $(\al_{p,1},\al_{p,2},\beta_{p,1},\beta_{p,2})\in\R^4\setminus\{0\}$. 
As a consequence,
\be\lab{3-10-10} \la_{p,4}\le\max_\LR \f{\int_{B_1(0)}|\nabla f_1+\nabla f_2|^2\rd x}{p\int_{B_1(0)}|u_p|^{p-1}(f_1+f_2)^2\rd x}\ee
where
    \[\LR=\{(\al_{p,1},\al_{p,2},\beta_{p,1},\beta_{p,2})\in\R^4\setminus\{0\}:~\sum_{j=1}^2(\al_{p,j}^2+\beta_{p,j}^2)=1\}.\]
Using \eqref{3-10-1}-\eqref{3-10-3}, we obtain
    \[\int_{B_1(0)}|\nabla f_1|^2=\sum_{j=1}^2\al_{p,j}^2\int_{B_1(0)}|\nabla\phi_j|^2u_p^2+\phi_j^2|u_p|^{p+1}\rd x, \]
    \[\int_{B_1(0)}|\nabla f_2|^2=\sum_{j=1}^2\beta_{p,j}^2\int_{B_1(0)}|\nabla\phi_j|^2|\pa_{x_1}u_p|^2+p\phi_j^2|u_p|^{p-1}|\pa_{x_1}u_p|^2\rd x, \]
and
    \[\int_{B_1(0)}\nabla f_1\cdot\nabla f_2=\sum_{j=1}^2\al_{p,j}\beta_{p,j}\int_{B_1(0)}|\nabla\phi_j|^2u_p\pa_{x_1}u_p+\f{p+1}{2}\phi_j^2|u_p|^{p-1}u_p\pa_{x_1}u_p\rd x.\] 
Therefore, we obtain
\be\lab{3-10-11}\begin{aligned}  
    &\int_{B_1(0)}|\nabla f_1+\nabla f_2|^2\rd x \\
    &=\sum_{j=1}^2\al_{p,j}^2\int_{B_1(0)}\phi_j^2|u_p|^{p+1}\rd x +(p+1)\sum_{j=1}^2\al_{p,j}\beta_{p,j}\int_{B_1(0)}\phi_j^2|u_p|^{p-1}u_p\pa_{x_1}u_p\rd x \\
    &\quad +p\sum_{j=1}^2\beta_{p,j}^2\int_{B_1(0)}\phi_j^2|u_p|^{p-1}|\pa_{x_1}u_p|^2\rd x +\sum_{j=1}^2\int_{B_1(0)}|\nabla\phi_j|^2(\al_{p,j}u_p+\beta_{p,j}\pa_{x_1}u_p)^2\rd x \\
    &=:I_1+(p+1)I_2+pI_3+I_4.
\end{aligned}\ee 
Furthermore, by direct computations, we have
\be\lab{3-10-12}  \int_{B_1(0)}|u_p|^{p-1}(f_1+f_2)^2\rd x=I_1+2I_2+I_3. \ee 
Hence from \eqref{3-10-10}-\eqref{3-10-12}, we obtain
\be\lab{3-10-13} \la_{p,4}\le\max_\LR\f{I_1+(p+1)I_2+pI_3+I_4}{p(I_1+2I_2+I_3)}=1+\max_\LR\f{(1-p)I_1+(1-p)I_2+I_4}{pI_1+2pI_2+pI_3}. \ee 

Now, we estimate each $I_j$. By \eqref{est-wpi-3} and \eqref{out-u}, we obtain
\be\lab{3-10-14}\begin{aligned}
    I_1 &=\sum_{j=1}^2\al_{p,j}^2\sbr{\int_{B_r(x_{p,j})}|u_p|^{p+1}\rd x+O\sbr{\f{C^{p+1}}{p^{p+1}}}}\\
    &=\sum_{j=1}^2\al_{p,j}^2\sbr{\f{|u_p(x_{p,j})|^2}{p}\int_{B_{\f{r}{\ve_{p,j}}}(0)}\abs{1+\f{w_{p,j}}{p}}^{p+1}\rd y+O\sbr{\f{C^{p+1}}{p^{p+1}}}}\\
    &=\f{e}{p}\sum_{j=1}^2\al_{p,j}^2\sbr{\int_{\R^2}e^U\rd y+o(1)}\\
    &=\f{8\pi e}{p}\sum_{j=1}^2\al_{p,j}^2(1+o(1)),
\end{aligned}\ee
and similarly,
\be\lab{3-10-15}\begin{aligned}
    I_2 &=\sum_{j=1}^2\al_{p,j}\beta_{p,j}\sbr{\int_{B_r(x_{p,j})}\pa_{x_1}\sbr{\f{|u_p|^{p+1}}{p+1}}\rd x+O\sbr{\f{C^{p+1}}{p^{p+1}}}}\\
    &=\f{1}{p+1}\sum_{j=1}^2\al_{p,j}\beta_{p,j}\sbr{\int_{\pa B_r(x_{p,j})}|u_p|^{p+1}\nu_1\rd \sigma+O\sbr{\f{C^{p+1}}{p^{p+1}}}}\\
    &=O\sbr{\f{C^{p+1}}{p^{p+1}}},
\end{aligned}\ee
\be\lab{3-10-16}\begin{aligned}
    I_3 &=\sum_{j=1}^2\beta_{p,j}^2\sbr{\int_{B_r(x_{p,j})}|u_p|^{p-1}|\pa_{x_1}u_p|^2\rd x+O\sbr{\f{C^{p+1}}{p^{p+1}}}}\\
    &=\sum_{j=1}^2\beta_{p,j}^2\sbr{\f{|u_p(x_{p,j})|^2}{p^3\ve_{p,j}^2}\int_{B_{\f{r}{\ve_{p,j}}}(0)}\abs{1+\f{w_{p,j}}{p}}^{p-1}|\pa_{y_1}w_{p,j}|^2\rd y+O\sbr{\f{C^{p+1}}{p^{p+1}}}}\\
    &\ge\sum_{j=1}^2\f{e\beta_{p,j}^2}{p^3\ve_{p,j}^2}\sbr{\int_{\R^2}e^U|\pa_{y_1}U|^2\rd y+o(1)}+O\sbr{\f{C^{p+1}}{p^{p+1}}},
\end{aligned}\ee
and
\be\lab{3-10-17} 0\le I_4\le C\sum_{j=1}^2\int_{B_{3r}(X_j)\setminus B_{2r}(X_j)}(\al_{p,j}u_p+\beta_{p,j}\pa_{x_1}u_p)^2\rd x=O\sbr{\f{1}{p^2}}.\ee
Write $A_p=\sum_{j=1}^2\al_{p,j}^2$ and $B_p=\sum_{j=1}^2\beta_{p,j}^2$, so that $A_p+B_p=1$. All estimates in \eqref{3-10-14}--\eqref{3-10-17} are uniform on $\LR$. They give constants $c,C>0$, independent of the coefficients, such that
\[(1-p)I_1+(1-p)I_2+I_4\le-cA_p+\frac C{p^2},\]
and
\[p(I_1+2I_2+I_3)\ge c\left(A_p+\frac{B_p}{p^2\ve_p^2}\right).\]
Here the superexponentially small term $I_2$ is absorbed using $A_p+B_p=1$ and the comparability of $\ve_{p,1},\ve_{p,2}$ and $\ve_p$ from (P4). If $A_p\ge 2C/(cp^2)$, the first expression is nonpositive. Otherwise $B_p\ge1/2$ for sufficiently large $p$, and the quotient of these two expressions is at most $C\ve_p^2$. Substitution into \eqref{3-10-13} proves \eqref{la_34} uniformly over $\LR$.
\end{proof}

\begin{lemma}\lab{lemma-3-11-2}
For $l=3,4$ and $j=1,2$, we have $b_{l,j}=0$ and
\be\lab{in-vplj-1}\ti v_{p,l,j}=a_{l,j}\f{x_1}{8+|x|^2}+o(1)\quad\text{in}~C_\loc^1(\R^2),~\text{as}~p\to\iy,\ee
where $a_{l,j},b_{l,j}$ are constants in \eqref{in-vplj}. Moreover, we have
\be\lab{3-11-1} \begin{aligned}
    \boldsymbol{a}_l=(a_{l,1},a_{l,2})\neq0,~\forall~l=3,4, \quad\text{and}\quad   \boldsymbol{a}_3\cdot \boldsymbol{a}_4=0.
\end{aligned}\ee
\end{lemma}
\begin{proof}
Since $\la_{\iy,l}=1$, we have \eqref{in-vplj}. If $b_{l,j}\neq0$ for some $j\in\{1,2\}$, Lemma~\ref{lemma-4} gives $\la_{p,l}=1+\f{6}{p}+o\sbr{\f1p}$, which contradicts \eqref{la_34}. Thus $b_{l,j}=0$ for $j=1,2$, and \eqref{in-vplj-1} holds. Lemma~\ref{lemma-2} then gives \eqref{3-11-1}.
\end{proof}

The local Pohozaev identities determine the first nonzero shifts of $\la_{p,3}$ and $\la_{p,4}$.
\begin{lemma}\lab{lemma-}
The eigenvalues satisfy
\be\lab{3-11-18} \begin{aligned}
        \la_{p,3}&=1+6(3+\sqrt 5)e^{-2\pi\Psi(\bX)}\ve_p^2+o(\ve_p^2),\\
        \la_{p,4}&=1+12(5+2\sqrt5)e^{-2\pi\Psi(\bX)}\ve_p^2+o(\ve_p^2).
\end{aligned}\ee 
\end{lemma}
\begin{proof}
Fix a sufficiently small $d>0$. Rescaling \eqref{B_plj} and using \eqref{est-wpi-3} and \eqref{in-vplj-1}, we obtain
\be\lab{3-11-2} \begin{aligned} 
    B_{p,l,j}&=\f{\ve_{p,j}}{p}\int_{B_{\f{d}{\ve_{p,j}}}(0)}y_1\abs{1+\f{w_{p,j}}{p}}^{p-1}\ti v_{p,l,j}\rd y\\
    &=\f{\ve_{p,j}}{p}\sbr{a_{l,j}\int_{B_{\f{d}{\ve_{p,j}}}(0)}\f{y_1^2}{8+|y|^2}e^{U(y)}\rd y+o(1)}\\
    &=2\pi a_{l,j}\f{\ve_{p,j}}{p}+o\sbr{\f{\ve_{p}}{p}}.
\end{aligned} \ee
Combining \eqref{out-vpl-1} and \eqref{3-11-2}, we obtain
\be\lab{3-11-3} v_{p,l}(x)=\la_{p,l}p \sum_{j=1}^2 A_{p,l,j}G(x_{p,j},x)+2\pi\la_{p,l}\sum_{j=1}^2 a_{l,j}\ve_{p,j}\pa_1G(x_{p,j},x)+o(\ve_p), \ee
in $C^1(\ov{B_1(0)}\setminus(B_{2d}(X_1)\cup B_{2d}(X_2)))$.

First, we use the Pohozaev identity \eqref{Pj-1} to estimate $A_{p,l,j}$.
By Lemma \ref{lemma-3-11-2} we have $b_{l,j}=0$, and hence \eqref{3-11-5} implies $p^2A_{p,l,j}=o(1)$ for each $j=1,2$. 
Using \eqref{Pj}, \eqref{out-u} and \eqref{3-11-3}, we obtain
\[\begin{aligned}
    P_j(u_p,v_{p,l})&=p\la_{p,l}\sum_{m,k=1}^2C_{p,m}A_{p,l,k}P_j(G(x_{p,m},\cdot),G(x_{p,k},\cdot))\\
    &\quad+2\pi\la_{p,l}\sum_{m,k=1}^2C_{p,m}a_{l,k}\ve_{p,k}P_j(G(x_{p,m},\cdot),\pa_1G(x_{p,k},\cdot))+o\sbr{\f{\ve_p}{p}}.
\end{aligned}\]
Lemma \ref{lemma-3-11-1} and \eqref{Cpi} give
\be\lab{3-11-6}\begin{aligned}
    P_j(u_p,v_{p,l})&=-\f{1}{2\pi}\la_{p,l}pC_{p,j}A_{p,l,j}+o\sbr{\f{\ve_p}{p}}\\
    &\quad +2\pi\la_{p,l}\sbr{\f12C_{p,j}a_{l,j}\ve_{p,j}\pa_{x_1}R(x_{p,j})+C_{p,j'}a_{l,j}\ve_{p,j}D_1G(x_{p,j'},x_{p,j})}\\
    &=-\f{1}{2\pi}\la_{p,l}pC_{p,j}A_{p,l,j}+(-1)^j\f{8\pi^2\sqrt e}{p}\la_{p,l}a_{l,j}\ve_{p,j}\pa_{a_{2j-1}}\Psi(\bX)+o\sbr{\f{\ve_p}{p}}
\end{aligned}\ee
where $j'=3-j$ and $\bX=(X_1,X_2)\in\R^4$. Here $a_{2j-1}$ is the first coordinate of the $j$th point in the four real coordinates of $\Psi$. Since $\bX$ is a critical point of $\Psi$, we obtain
\be\lab{3-11-7} P_j(u_p,v_{p,l})=(-1)^j(4\sqrt e+o(1))A_{p,l,j}+o\sbr{\f{\ve_p}{p}}. \ee 
On the other hand, by \eqref{out-u} and \eqref{3-11-3}, we obtain
\be\lab{3-11-8} \int_{\pa B_d(x_{p,j})}|u_p|^{p-1}u_pv_{p,l} (x-x_{p,j})\cdot\nu\rd\sigma=O\sbr{\f{\ve_p}{p^2}+\f1p(|A_{p,l,1}|+|A_{p,l,2}|)}. \ee
The equations for $u_p$ and $v_{p,l}$ give
\be\lab{3-11-9}\begin{aligned} 
    (p\la_{p,l}-1)\int_{B_d(x_{p,j})}|u_p|^{p-1}u_pv_{p,l}\rd x
    &=\int_{B_d(x_{p,j})}\Delta u_p v_{p,l}-\Delta v_{p,l}u_p\rd x\\
    &=\int_{\pa B_d(x_{p,j})}v_{p,l}\pa_\nu u_p -u_p\pa_\nu v_{p,l}\rd \sigma,
\end{aligned}\ee
Using \eqref{out-u} and \eqref{3-11-3} on the boundary, we obtain
\be\lab{3-11-10} \int_{B_d(x_{p,j})}|u_p|^{p-1}u_pv_{p,l}\rd x=O\sbr{\f{\ve_p}{p^2}+\f1p(|A_{p,l,1}|+|A_{p,l,2}|)}. \ee 
Using the rescaling in \eqref{ti-v}, the integral below is $u_p(x_{p,j})/p$ times an integral converging to $\int e^UV_{l,j}y\cdot\nabla U=0$. The limit vanishes by oddness because $b_{l,j}=0$. Thus
\be\lab{3-11-11}  p\int_{B_d(x_{p,j})}|u_p|^{p-1}v_{p,l}(x-x_{p,j})\cdot\nabla u_p\rd x=o\sbr{\f1p}. \ee
Substituting \eqref{3-11-7}, \eqref{3-11-8}, \eqref{3-11-10} and \eqref{3-11-11} into \eqref{Pj-1} and taking absolute values for each $j$, we sum over $j=1,2$ and absorb the term $C\sum_j|A_{p,l,j}|/p$ into the left-hand side for large $p$. This gives
\be\lab{3-11-12} |A_{p,l,j}|=o\sbr{\f{\ve_p+|\la_{p,l}-1|}{p}},\quad\forall~j=1,2. \ee 

Next we use the Pohozaev identity \eqref{Qj-1} to estimate $\la_{p,l}-1$.
Expanding $Q_j$ as in \eqref{3-11-6} and using \eqref{3-11-12}, we obtain
\be\lab{3-11-13} \begin{aligned}
    Q_j(u_p,v_{p,l}) &= 2\pi\la_{p,l}\sum_{m,k=1}^2C_{p,m}a_{l,k}\ve_{p,k}Q_j\big(G(x_{p,m},\cdot),\pa_1 G(x_{p,k},\cdot)\big)+o\sbr{\f{\ve_p+|\la_{p,l}-1|}{p}} \\
    &=(-1)^j\f{8\pi^2\sqrt e}{p}\la_{p,l}a_{l,j}\ve_{p,j}\sbr{-\pa^2_{11}R(x_{p,j})+2D_{11}^2G(x_{p,j'},x_{p,j})} \\
        &\quad +(-1)^{j'}\f{8\pi^2\sqrt e}{p}\la_{p,l}a_{l,j'}\ve_{p,j'} 2D_1\pa_1G(x_{p,j'},x_{p,j}) +o\sbr{\f{\ve_p+|\la_{p,l}-1|}{p}} \\
    &=\f{8\pi^2\sqrt e}{p}\la_{p,l}\sum_{k=1}^2(-1)^ka_{l,k}\ve_{p,k}\pa_{a_{2j-1}}\pa_{a_{2k-1}}\Psi(\bX)+o\sbr{\f{\ve_p+|\la_{p,l}-1|}{p}}. 
\end{aligned}\ee
where $j'\in\{1,2\}$ and $j'\neq j$.
Using \eqref{out-u}, \eqref{3-11-3} and \eqref{3-11-12}, we obtain
\be\lab{3-11-14} \int_{\pa B_d(x_{p,j})}|u_p|^{p-1}u_pv_{p,l} \nu_1\rd\sigma =o\sbr{\f{\ve_p+|\la_{p,l}-1|}{p}}. \ee
The same rescaling, now using $\int e^UV_{l,j}\partial_1U=-\pi a_{l,j}/3$, gives
\be\lab{3-11-15}  p\int_{B_d(x_{p,j})}|u_p|^{p-1}v_{p,l}\pa_{x_1}u_p\rd x=\f{(-1)^j}{p\ve_{p,j}}\sbr{\f{\pi\sqrt e}{3}a_{l,j}+o(1)}. \ee
Thus, substituting \eqref{3-11-13}-\eqref{3-11-15} into \eqref{Qj-1}, we obtain
\be\lab{3-11-16} \begin{aligned} 
    (\la_{p,l}-1)\sbr{\f{\pi\sqrt e}{3}a_{l,j}+o(1)} 
    &=8\pi^2\sqrt e\sum_{k=1}^2(-1)^{k-j}a_{l,k}\ve_{p,j}\ve_{p,k}\pa_{a_{2j-1}}\pa_{a_{2k-1}}\Psi(\bX)\\
    &\quad +o\sbr{\ve_p^2+\ve_p|\la_{p,l}-1|},\quad\forall~j=1,2,
\end{aligned}\ee 
which implies 
\be \lab{3-11-17} |\la_{p,l}-1|=O(\ve_p^2), \ee 
because at least one component of $\boldsymbol a_l$ is nonzero by \eqref{3-11-1}.

By \eqref{3-11-16} and \eqref{3-11-17}, we obtain
    \[a_{l,j}(\la_{p,l}-1)=24\pi\sum_{k=1}^2(-1)^{k-j}a_{l,k}\ve_{p,j}\ve_{p,k}\pa_{a_{2j-1}}\pa_{a_{2k-1}}\Psi(\bX)+o(\ve_p^2), \]
for any $l=3,4$ and $j=1,2$.
In other words, we have
    \[(\la_{p,l}-1)\boldsymbol{a}_l=24\pi E_pSTSE_p\boldsymbol{a}_l+o(\ve_p^2),\quad \forall~l=3,4,\]
where $E_p=\mathrm{diag}(\ve_{p,1},\ve_{p,2})$, $S=\mathrm{diag}(1,-1)$ and
\[T=\sbr{\begin{matrix} \pa_{a_1a_1}^2\Psi(\bX) & \pa_{a_1a_3}^2\Psi(\bX) \\ \pa_{a_3a_1}^2\Psi(\bX) & \pa_{a_3a_3}^2\Psi(\bX) \end{matrix} }
    =\f{1}{8\pi}\sbr{ \begin{matrix}
    13+5\sqrt5 & -7-3\sqrt5 \\
    -7-3\sqrt5 & 13+5\sqrt5
    \end{matrix} }.\]
At the symmetric critical pair $\bX$, one has $\Psi_1(\bX)=\Psi_2(\bX)=\Psi(\bX)/2$, and the scale relations \eqref{ve_p} and \eqref{ve_p-1} imply $\ve_{p,j}=e^{-\pi\Psi(\bX)}\ve_p+o(\ve_p)$ for $j=1,2$. Hence
    \[\f{24\pi E_pSTSE_p}{\ve_p^2}\to24\pi e^{-2\pi\Psi(\bX)}STS.\]
The two eigenvalues of the limit matrix are
    \[6(3+\sqrt 5)e^{-2\pi\Psi(\bX)},\quad 12(5+2\sqrt5)e^{-2\pi\Psi(\bX)}.\]
Using the global orthogonality \eqref{3-11-1} and the ordering $\la_{p,3}\le \la_{p,4}$, we obtain the refined expansions
    \[\begin{aligned}
        \la_{p,3}&=1+6(3+\sqrt 5)e^{-2\pi\Psi(\bX)}\ve_p^2+o(\ve_p^2),\\
        \la_{p,4}&=1+12(5+2\sqrt5)e^{-2\pi\Psi(\bX)}\ve_p^2+o(\ve_p^2).
    \end{aligned}\]
This completes the proof.
\end{proof}

We now prove Proposition \ref{thm-index}.
\begin{proof}[Proof of Proposition \ref{thm-index}]
Along the arbitrary sequence fixed at the beginning of this section, \eqref{3-10-4} and \eqref{3-11-18} give
\[0<\la_{p,1}<\la_{p,2}<1<\la_{p,3}<\la_{p,4}\]
for sufficiently large $p$. Since the remaining eigenvalues are no smaller than $\la_{p,4}$, the quadratic form has exactly two negative directions and no kernel in $\HR_e$. The orientation change described after \eqref{assumption-1} preserves the spectrum. If the proposition failed for arbitrarily large $p$, a sequence of such counterexamples would contradict this conclusion.
\end{proof}

\vskip0.2in
\section{Leray--Schauder degree in the symmetric space}\lab{sec-degree}

Esposito, Musso and Pistoia \cite{Esposito-Musso-Pistoia-2007} used finite-dimensional reduction to construct a family of low-energy nodal solutions of \eqref{equ-1}. We show that every even low-energy nodal solution admits a two-peak representation of the same form. This representation reduces the computation of the Leray--Schauder degree to a finite-dimensional problem.

\begin{proposition}\lab{thm-degree}
For all sufficiently large $p$,
\[
 \deg(\mathrm{id}-\mathcal T_p,\DR_p\cap\mathcal H_e,0)=2,
\]
where $\TR_p$ and $\DR_p$ are given by \eqref{TR} and \eqref{define-d_p}.
\end{proposition}

Blow-up analysis has been used to study the Leray--Schauder degree for mean-field equations \cite{Li-1999,Chen-Lin-2003}.
We adapt the finite-dimensional reduction and degree argument of \cite[Sections 5--7]{Ianni-Luo-Yan-2025} to the signed two-peak approximation in the even subspace.
By (P3) of Section~\ref{sec-2}, the signed profile $(-1)^{j-1}u_p$ satisfies the positive Lane--Emden equation in a fixed neighborhood of the $j$th peak. We give the computations needed for opposite peak signs and for the restriction to the even subspace. Related estimates appear in \cite{Esposito-Musso-Pistoia-2007}.

For an even low-energy nodal solution of \eqref{equ-1}, the peak locations $x_{p,1}$ and $x_{p,2}$ both lie on the $x_1$-axis. Hence
    \[\text{either}~(x_{\iy,1},x_{\iy,2})=(X_1,X_2),\quad \text{or}~(x_{\iy,1},x_{\iy,2})=-(X_1,X_2),\]
where $X_1,X_2$ are given by \eqref{X}.
Set $\boldsymbol X=(X_1,X_2)$. For $a=(a_1,0)$ and $r>0$, write
\[
 L_r(a)=\{(t,0)\in\mathbb R^2:|t-a_1|<r\}.
\]
Choose $d_0>0$ so small that the closed intervals $\overline{L_{2d_0}(X_1)}$ and $\overline{L_{2d_0}(X_2)}$ are disjoint and contained in $B_1(0)$. 

The radial functions $w_0$ and $w_1$ from Section~\ref{sec2.2} satisfy the following expansions as $r\to\infty$ (see \cite[(2.4)--(2.7)]{Esposito-Musso-Pistoia-2007}).
\be\lab{sec5-radial-expansion}
 w_i(r)=C_i\log r+\kappa_i+
 O\big(r^{-2}(1+\log^8r)\big),\qquad i=0,1.
\ee
Moreover, we have $C_0=12$.
The radial error equations \eqref{3-12-0} and \eqref{eq-w1-correct} imply the corresponding differentiated estimates of \eqref{sec5-radial-expansion}. 
For each $\bxi=(\xi_1,\xi_2)\in L_1(0)\times L_1(0)$ close to $\bX$ or $-\bX$, we obtain the scales $\bar\bmu=(\bar\mu_{p,1},\bar\mu_{p,2})\in(0,\infty)^2$ by solving the following system.
\be\lab{Fpj}\begin{split}
0=F_{p,j}(\bar{\boldsymbol\mu},\boldsymbol\xi)
={}&-\log(64\bar\mu_{p,j}^4)
 +8\pi A_pR(\xi_j)
 +\left(\frac{C_0}{p}+\frac{C_1}{p^2}\right)
       \log(\bar\mu_{p,j}e^{-p/4})\\
 &-8\pi A_p
   \left(\frac{\bar\mu_{p,j}}{\bar\mu_{p,j'}}\right)^{2/(p-1)}G(\xi_j,\xi_{j'})
 -\frac{\kappa_0}{p}-\frac{\kappa_1}{p^2},
\end{split}\ee
where $j'\in\{1,2\}$ with $j'\neq j$ and $A_p=1-\frac{C_0}{4p}-\frac{C_1}{4p^2}$.
We point out that the definition of $F_{p,j}(\bar\bmu,\bxi)$ is inspired by \cite[(2.11)]{Esposito-Musso-Pistoia-2007}. 
In logarithmic variables $m_j=\log\bar\mu_{p,j}$, the derivative of this system is $-4I+O(p^{-1})$, uniformly on bounded sets of $m_j$. Its limit has the unique solution
\[ m_j=-2\pi\Psi_j(\boldsymbol\xi)-\frac32\log2-\frac34.\]
The implicit function theorem gives a unique solution near this limit, depending smoothly on $\boldsymbol\xi$, and
\be\lab{3-12-20}
 \bar\mu_{p,j}=\f{1}{2\sqrt 2}e^{-2\pi\Psi_j(\boldsymbol\xi)-3/4}\big(1+O(p^{-1})\big).
\ee
In particular, for any fixed $d_0>0$ there exists $c>0$ such that $\f1c\le\bar\mu_{p,j}\le c$ for $\bxi\in L_{d_0}(X_1)\times L_{d_0}(X_2)$ or $\bxi\in L_{d_0}(-X_1)\times L_{d_0}(-X_2)$. 

The exact algebraic form of \eqref{Fpj} is
\[\begin{aligned}
\log\bar\mu_{p,j}&=2\pi\left[R(\xi_j)-\left(\frac{\bar\mu_{p,j}}{\bar\mu_{p,j'}}\right)^{2/(p-1)}G(\xi_j,\xi_{j'})\right]\\
&\quad-\frac{2\log8+C_0/4+C_1/(4p)+\kappa_0/p+\kappa_1/p^2}{4A_p}.
\end{aligned}\]

For $j=1,2$, we set
\be\lab{3-12-21}  \bar\ve_{p,j}=\bar\mu_{p,j}e^{-\f{p}{4}},\qquad a_{p,j}=\frac1{p^{p/(p-1)}\bar\varepsilon_{p,j}^{2/(p-1)}}.\ee 
and inspired by \cite{Esposito-Musso-Pistoia-2007}, we define 
\be\lab{Wpj}
 \bar W_{p,j}(x)=a_{p,j}\left[p+U\left(\frac{x-\xi_j}{\bar\varepsilon_{p,j}}\right)+\frac1p w_0\left(\frac{x-\xi_j}{\bar\varepsilon_{p,j}}\right)+\frac1{p^2}w_1\left(\frac{x-\xi_j}{\bar\varepsilon_{p,j}}\right)\right].
\ee
It is easy to check that $W_{p,j}$ is symmetric since $\xi_j\in L_{d_0}(\pm X_j)$, and hence $PW_{p,j}$ is symmetric, i.e., $PW_{p,j}\in\HR_e$.
For a smooth function $v$ on $\overline{B_1(0)}$, its projection $Pv$ is the unique zero-boundary solution of $\Delta Pv=\Delta v$. The expansion on $\partial B_1(0)$ and the maximum principle give
\be\lab{sec5-projection}
    P\bar W_{p,j}(x)=\bar W_{p,j}(x)-a_{p,j}K_{p,j}+8\pi a_{p,j}A_pH(x,\xi_j)+O(a_{p,j}p^N\bar\varepsilon_{p,j}^2),
\ee
where $N$ is a fixed integer and
    \[K_{p,j}=p+\log64+4A_p\log\bar\varepsilon_{p,j}+\frac{\kappa_0}{p}+\frac{\kappa_1}{p^2}.\]
The error is exponentially small in $p$. Since
$a_{p,j'}/a_{p,j}=(\bar\mu_{p,j}/\bar\mu_{p,j'})^{2/(p-1)}$,
equations \eqref{Fpj} are precisely
    \[ K_{p,j}=8\pi A_p\left[R(\xi_j)-\frac{a_{p,j'}}{a_{p,j}}G(\xi_j,\xi_{j'})\right].\]
Thus the constants $\kappa_i$ disappear from the projected sum after this matching condition is imposed.

Let $s_j=(-1)^{j-1}$ and
\be\lab{Zpk-Wap}
    Z_{p,j}=s_jP\bar W_{p,j},\qquad
    W_{\boldsymbol\alpha,p}=\sum_{j=1}^2\alpha_j Z_{p,j}.
\ee
It is worth noting that, in contrast to the positive solution case, we adopt the basis $(-1)^{j-1}P\ov W_{p,j}$ for the nodal case.
At $\boldsymbol\alpha=(1,1)$, \eqref{sec5-projection} and \eqref{Fpj} show that near $\xi_j$,
\[
 s_jW_{\boldsymbol1,p}(x)=\bar W_{p,j}(x)
 +O\big(a_{p,j}|x-\xi_j|\big)+O(p^Ne^{-p/2}).
\]
On every fixed closed set separated from the centers $\xi_j$, they give
\be\lab{sec5-outer-W}
 W_{\boldsymbol\alpha,p}=8\pi A_p\sum_{j=1}^2\alpha_js_ja_{p,j}G(\cdot,\xi_j)+O(p^Ne^{-p/2})
\ee
in $C^1$. In particular, $a_{p,j}=\sqrt e\,p^{-1}(1+O(p^{-1}\log p))$.
For each $\boldsymbol\xi$, let
\be\lab{ORe}
 \mathcal O_e(\boldsymbol\xi)=\left\{w\in\mathcal H_e:\langle Z_{p,j},w\rangle=\langle\partial_{\xi_{j,1}}Z_{p,j},w\rangle=0,\ j=1,2\right\}.
\ee
Write $Y_{p,j}=\partial_{\xi_{j,1}}Z_{p,j}$. The computations in \cite[(6.8)--(6.10) and (7.4)--(7.7)]{Ianni-Luo-Yan-2025}, applied to \eqref{sec5-projection}, give, in our normalization,
\[ \|Z_{p,j}\|^2=\frac{8\pi e}{p}+o(p^{-1}),\quad \|Y_{p,j}\|^2=\frac{e+o(1)}{p^2\bar\varepsilon_{p,j}^2}\int_{\mathbb R^2}|\nabla\partial_1U|^2\,dy,\]
and
\[
\langle Z_{p,1},Z_{p,2}\rangle=O(p^{-2}),\quad \frac{\langle Z_{p,i},Y_{p,j}\rangle}{\|Z_{p,i}\|\|Y_{p,j}\|}=o(1), \quad \frac{\langle Y_{p,1},Y_{p,2}\rangle}{\|Y_{p,1}\|\|Y_{p,2}\|}=o(1).
\]

To compute the degree, we work in $\mathcal X_e=\mathcal H_e\cap C_0(\overline{B_1(0)})$ with the norm $\|w\|+\|w\|_\infty$.
The compactness of $\mathcal T_p\colon\mathcal H_e\to\mathcal X_e$ for fixed $p$ follows from the Sobolev embeddings in dimension two, the Dirichlet $W^{2,q}$ estimate with $q>2$, and the compact embedding of $W^{2,q}$ into $C^1$. 
For fixed $p$, every solution in $\overline{\mathcal E_p}$ satisfies the PDE identity
\[\|u\|^2=\int_{B_1(0)}|u|^{p+1}\,dx\le\frac{16\pi e}{p}.\]
Thus, all these solutions belong to the bounded set $\DR_p$. Elliptic estimates give an $L^\infty$ bound for these solutions. Choose $R_p$ strictly larger than that bound and set
    \[\mathcal U_p=\DR_p\cap\mathcal X_e\cap\{u:\|u\|_\infty<R_p\}.\]
This set is bounded in $\mathcal X_e$, and its additional boundary contains no zero of $I-\TR_p$. Let $j\colon\mathcal X_e\to\mathcal H_e$ be the inclusion.
The two realizations of the compact operator are $j\circ\TR_p$ and $\TR_p\circ j$. The commutativity property of Leray--Schauder degree gives
\be\lab{sec5-degree-spaces}
 \deg_{\mathcal H_e}(I-j\circ\TR_p,\DR_p\cap\mathcal H_e,0)=\deg_{\mathcal X_e}(I-\TR_p\circ j,\mathcal U_p,0).
\ee

Fix $0<\tau_0<1/4$. Let $\mathcal S_p^\pm$ be defined by
\be \SR_p^+=\lbr{(\bal,\bxi,w):
\begin{aligned}
&\bal\in\R^2,\bxi\in L_1(0)^2,w\in\OR_e\cap C_0(\overline{B_1(0)}),\\
&|\al_j-1|\le\f{1}{p^{\f52}},|\xi_j-X_j|\le d_0,\forall~j=1,2,\\
&\nm{w}\le\f{1}{p^{2+\tau_0}},\nm{w}_{\iy}\le\f{1}{p^{2+\tau_0}}
\end{aligned}  }, \ee
and
\be \SR_p^-=\lbr{(\bal,\bxi,w):
\begin{aligned}
&\bal\in\R^2,\bxi\in L_1(0)^2,w\in\OR_e\cap C_0(\overline{B_1(0)}),\\
&|\al_j-1|\le\f{1}{p^{\f52}},|\xi_j+X_j|\le d_0,\forall~j=1,2,\\
&\nm{w}\le\f{1}{p^{2+\tau_0}},\nm{w}_{\iy}\le\f{1}{p^{2+\tau_0}}
\end{aligned}  }. \ee
The sets $\SR_p^{\pm}$ are parameter neighborhoods of the configurations $\pm\bX$.
Let $\eta_j(x)=\eta(x-\xi_j)$, where $\eta$ is a smooth radial cutoff equal to one on $B_r(0)$ and zero outside $B_{r+d}(0)$.
For $u\in\mathcal X_e$, we define the stress functional by
\be\lab{sec5-stress}
\mathcal P_j(u)=\int_{B_1(0)}\left[\partial_1u\,\nabla u\cdot\nabla\eta_j-\frac12|\nabla u|^2\partial_1\eta_j+\frac{|u|^{p+1}}{p+1}\partial_1\eta_j\right]dx.
\ee
This is continuous in $\mathcal X_e$. For smooth $u$, integration by parts gives
    \[\mathcal P_j(u)=\langle J_p'(u),\eta_j\partial_1u\rangle.\]
Then, we define the following maps on $(\bal,\bxi,w)\in\SR_p^+\cup\SR_p^-$:
\[\begin{aligned}
 \mathcal A_{p,1,j}&=\langle J_p'(W_{\boldsymbol\alpha,p}+w),Z_{p,j}\rangle,\\
 \mathcal A_{p,2,j}&=\mathcal P_j(W_{\boldsymbol\alpha,p}+w),\\
 \mathcal A_{p,3}&=\Pi_{\boldsymbol\xi}
       \big[W_{\boldsymbol\alpha,p}+w-\mathcal T_p(W_{\boldsymbol\alpha,p}+w)\big].
\end{aligned}\]
and $\AR_p=(\AR_{p,1,1},\AR_{p,1,2},\AR_{p,2,1},\AR_{p,2,2},\AR_{p,3})$. Then we have

\begin{lemma}\lab{lemma-3-12-1}
For all sufficiently large $p$, every solution $u_p\in\ER_p\cap\HR_e\cap\CR_p$ admits a representation $u_p=W_{\bal,p}+w$ with $(\bal,\bxi,w)\in\SR_p^+\cup\SR_p^-$ satisfying all bounds strictly and $\AR_p(\bal,\bxi,w)=0$.
Moreover, the representation is unique in a sufficiently small coordinate neighborhood, and
\be\lab{sec5-covering}
 |\alpha_j-1|=o(p^{-5/2}),\quad|\xi_j\mp X_j|=o(1),\quad\|w\|+\|w\|_\infty=o(p^{-2-\tau_0}).
\ee
Conversely, a zero of $\mathcal A_p$ in these parameter sets gives a solution in $\mathcal E_p$.
\end{lemma}
\begin{proof}
We adapt the argument of \cite[Proposition 7.1]{Ianni-Luo-Yan-2025} to signed peaks.
Start with the true centers $x_{p,j}$, scales $\varepsilon_{p,j}$, and $\mu_{p,j}=\varepsilon_{p,j}e^{p/4}$. 
In this paragraph a superscript $*$ denotes the unprojected profile \eqref{Wpj} using these true scales, before imposing \eqref{Fpj}. 
We write the source difference as $q_p$:
    \[-\Delta\sum_js_jP\bar W_{p,j}^*=|u_p|^{p-1}u_p-q_p. \]
Then, the weighted expansion \eqref{est-wpi-4} and the equations for $w_0,w_1$ imply that $q_p$ satisfies
\be\lab{sec5-source-true}
 |q_p(x_{p,j}+\varepsilon_{p,j}y)|
 \le\frac{C}{p^4\varepsilon_{p,j}^2}(1+|y|)^{-4+\delta},
 \quad\text{for}~ |y|\le r/\varepsilon_{p,j},
\ee
and
\be \int_{B_1\setminus\cup_jB_r(x_{p,j})}|q_p|\le C p^N e^{-p/2}, \ee
where $0<\delta<1$ is fixed after choosing the exponent in \eqref{est-wpi-4} sufficiently small. 
By applying Green's representation formula to $u_p-\sum_js_jP\bar W_{p,j}^*$, and using \eqref{sec5-source-true}, we get
\be\lab{sec5-true-global}
 \left\|u_p-\sum_js_jP\bar W_{p,j}^*\right\|_\infty\le Cp^{-3},\qquad \left\|u_p-\sum_js_jP\bar W_{p,j}^*\right\|\le Cp^{-7/2}.
\ee
Evaluation at $x_{p,j}$ in \eqref{sec5-projection}, using \eqref{sec5-true-global}, gives
$F_{p,j}(\boldsymbol\mu_p,\boldsymbol x_p)=O(p^{-2})$.
The uniformly invertible logarithmic Jacobian of $F_p$ implies
\[
 |\log\mu_{p,j}-\log\bar\mu_{p,j}(\boldsymbol x_p)|\le Cp^{-2}.
\]
Direct differentiation of the projected radial profiles gives
\[
\|\partial_{\log\mu_j}P\bar W_{p,j}\|+
\|\partial_{\log\mu_j}P\bar W_{p,j}\|_\infty\le Cp^{-1}.
\]
Thus replacement of the true scales by the matched scales leaves an error of order $p^{-3}$ in both norms.

Put $\varepsilon_j^0=\bar\varepsilon_{p,j}(\boldsymbol x_p)$ and define
\[\alpha_j=1+a_j,\qquad \xi_j=x_{p,j}+(p\varepsilon_j^0t_j,0).\]
For $q=(a_1,a_2,t_1,t_2)$, set
\[r(q)=u_p-W_{\boldsymbol\alpha,p},\qquad 
 \Phi_p(q)=\begin{pmatrix}
 (p\langle r(q),Z_{p,j}\rangle)_{j=1,2}\\
 (p\varepsilon_j^0\langle r(q),Y_{p,j}\rangle)_{j=1,2}
 \end{pmatrix}\in\mathbb R^4,\]
where $Y_{p,j}=\partial_{\xi_{j,1}}Z_{p,j}$ and all profiles and test functions are evaluated at $\boldsymbol\xi=\boldsymbol\xi(q)$. 
Thus $r(q)\in\mathcal O_e(\boldsymbol\xi)$ precisely when $\Phi_p(q)=0$. At $q=0$, the preceding remainder estimates give
\[\begin{aligned}
 |p\langle r(0),Z_{p,j}\rangle|
 &\le p\|r(0)\|_\infty\|\Delta Z_{p,j}\|_1\le Cp^{-3},\\
 |p\varepsilon_j^0\langle r(0),Y_{p,j}\rangle|
 &\le p\varepsilon_j^0\|r(0)\|_\infty\|\Delta Y_{p,j}\|_1
 \le Cp^{-3}.
\end{aligned}\]
The derivative satisfies
\[\begin{aligned}
 \partial_{\xi_{j,1}}W_{\boldsymbol1,p}&=Y_{p,j}+R_{p,j},\qquad
 \|R_{p,j}\|\le Cp^{-1}.
\end{aligned}\]
By \eqref{sec5-projection} and the oddness of $\partial_1U$, we get
\[p\|Z_{p,j}\|^2\to8\pi e,\qquad p^2(\varepsilon_j^0)^2\|Y_{p,j}\|^2\to e\int_{\mathbb R^2}|\nabla\partial_1U|^2\,dy,\]
and
\[p^2\varepsilon_j^0\big|\langle Z_{p,i},Y_{p,j}+R_{p,j}\rangle\big|\le C\sqrt p\,\varepsilon_j^0=o(1).\]
For $Q_{p,c}=\{q:|q|_\infty<c/p\}$, we have
\[
 \sup_{q\in Q_{p,c}}\|r(q)\|_\infty\le Cp^{-3}+\frac{Cc}{p}.
\]
Differentiating the profiles and the test functions therefore gives
\[\begin{aligned}
 \sup_{q\in Q_{p,c}}\|D\Phi_p(q)+M\|\le Cc+o(1)<\frac{\lambda_*}{2},
\end{aligned}\]
where $M=\operatorname{diag}(8\pi e I_2,e\int_{\mathbb R^2}|\nabla\partial_1U|^2\,dy I_2)$ and $\lambda_*=\min\{8\pi e,e\int_{\mathbb R^2}|\nabla\partial_1U|^2\,dy\}$,
after fixing $c>0$ small and then taking $p$ large. The map $q\mapsto q+M^{-1}\Phi_p(q)$ is a contraction of $\{|q|\le C_1p^{-3}\}$ into itself for a sufficiently large fixed $C_1$. Its fixed point satisfies $\Phi_p(q)=0$ and
\[
 |\alpha_j-1|\le Cp^{-3},\qquad
 |\xi_j-x_{p,j}|\le C\varepsilon_j^0p^{-2},\qquad
 \|r(q)\|+\|r(q)\|_\infty\le Cp^{-3}.
\]
This proves \eqref{sec5-covering}.

The signed profiles, (P5), and Lemma~\ref{lemma-3-6-2} give $x_{p,1}\in B_{2d_0}(\pm X_1)$ for a representation in $\mathcal S_p^\pm$. 
Compare the constructed representation $u_p=W_{\boldsymbol\alpha,p}+w$ with another representation $u_p=W_{\widetilde{\boldsymbol\alpha},p}+\widetilde w$ in the same parameter branch, and put
$h_{p,j}=(\widetilde\xi_j-\xi_j)/\bar\varepsilon_{p,j}$. Then
\[
 p\|W_{\boldsymbol\alpha,p}-W_{\widetilde{\boldsymbol\alpha},p}\|
 \le p(\|w\|+\|\widetilde w\|)\longrightarrow0.
\]
If $h_{p,j}\to h\ne0$ along a subsequence, rescaling in $B_{R\bar\varepsilon_{p,j}}(\xi_j)$ gives, for each fixed $R>0$,
\[\begin{aligned}
 0&=\lim_{p\to\infty}p^2
 \|W_{\boldsymbol\alpha,p}-W_{\widetilde{\boldsymbol\alpha},p}\|^2\ge e\int_{B_R}|\nabla U(y)-\nabla U(y-h)|^2\,dy>0,
\end{aligned}\]
which is a contradiction.
For $|h_{p,j}|\to\infty$, the second rescaled gradient tends to zero on $B_R$ and the same contradiction holds. Hence $h_{p,j}\to0$, so both parameter vectors belong to $Q_{p,c}$. On this convex set,
\[
 \frac{\lambda_*}{2}|q-\widetilde q|
 \le|\Phi_p(q)-\Phi_p(\widetilde q)|=0.
\]
Thus the representation is unique.

Conversely, suppose $\mathcal A_p(\boldsymbol\alpha,\boldsymbol\xi,w)=0$. Put $u_p=W_{\boldsymbol\alpha,p}+w$ and set
\[
\begin{aligned}
 B=(Z_{p,1},Z_{p,2},Y_{p,1},Y_{p,2}),\quad 
 T=(Z_{p,1},Z_{p,2},\eta_1\partial_1u_p,\eta_2\partial_1u_p).
\end{aligned}
\]
The equation $\mathcal A_{p,3}=0$ and elliptic regularity give
\[
 u_p-\mathcal T_p(u_p)=\sum_{k=1}^4\beta_kB_k.
\]
Testing with $T_i$ and using $\mathcal A_{p,1}=\mathcal A_{p,2}=0$, we obtain
\[
\begin{aligned}
 \sum_{k=1}^4\frac{\langle B_k,T_i\rangle}{\|B_i\|\|B_k\|}
 (\beta_k\|B_k\|)=0,\qquad
 \left(\frac{\langle B_k,T_i\rangle}{\|B_i\|\|B_k\|}\right)_{i,k=1}^4
 \to\operatorname{diag}(1,1,-1,-1).
\end{aligned}
\]
Here the terms containing $w$ satisfy
\[
\begin{aligned}
 &\|\partial_1(\eta_j\Delta Z_{p,k})\|_1\le\frac C{p\bar\varepsilon_{p,k}},\qquad\|\partial_1(\eta_j\Delta Y_{p,k})\|_1\le\frac C{p\bar\varepsilon_{p,k}^2},\\
 &\frac{|\langle B_k,\eta_j\partial_1w\rangle|}{\|B_k\|\|Y_{p,j}\|}\le\frac{\|w\|_\infty\|\partial_1(\eta_j\Delta B_k)\|_1}
 {\|B_k\|\|Y_{p,j}\|}=o(1).
\end{aligned}
\]
Therefore $\boldsymbol\beta=0$ and $u_p=\mathcal T_p(u_p)$. The peak profiles give $u_p(\xi_1)=\sqrt e+o(1)$ and $u_p(\xi_2)=-\sqrt e+o(1)$. By \eqref{sec5-projection},
\[
 \|Z_{p,j}\|^2=\frac{8\pi e}{p}+o(p^{-1}),\qquad
 \langle Z_{p,1},Z_{p,2}\rangle=O(p^{-2}).
\]
Since $\alpha_j\to1$ and $\|w\|\le p^{-2-\tau_0}$, we get
\[
 p\|u_p\|^2=p\|W_{\boldsymbol\alpha,p}\|^2+o(1)
 =16\pi e+o(1).
\]
The PDE identity therefore verifies the low-energy hypothesis before we invoke \eqref{eq-energy-window-expansion}, which gives
\[
 16\pi e-p\int_{B_1(0)}|u_p|^{p+1}\,dx
 =32\pi e\frac{\log p}{p}+O(p^{-1}).
\]
For sufficiently large $p$, this deficit lies strictly between zero and $48\pi e\log p/p$, so the solution belongs to $\mathcal E_p$.
\end{proof}

\begin{lemma}\lab{sec5-linear-lemma}
There exist $c,C>0$, independent of  $\boldsymbol\alpha,\boldsymbol\xi$ and $p$, such that
\be\lab{sec5-coercive}
 \|v\|^2-\int p|W_{\boldsymbol\alpha,p}|^{p-1}v^2\,dx\ge\frac c p\|v\|^2,
 \qquad \forall  v\in\mathcal O_e(\boldsymbol\xi),
\ee
and, uniformly for $t\in[0,1]$, if $v\in\mathcal O_e(\boldsymbol\xi)$ satisfies
\be\lab{sec5-inverse}
 -\Delta v-tp|W_{\boldsymbol\alpha,p}|^{p-1}v=h+\sum_j\left[b_j(-\Delta Z_{p,j})+c_j(-\Delta\partial_{\xi_{j,1}}Z_{p,j})\right],
\ee
then $\|v\|_\infty\le Cp\|h\|_*$, and the multipliers must satisfy
\be\lab{sec5-multiplier-bound}
 \sum_j\left(|b_j|+\frac{|c_j|}{\bar\varepsilon_{p,j}}\right)\le C\big(\|v\|_\infty+p\|h\|_*\big),
\ee
where the norm is defined by
\[\|h\|_* =\sup_{x\in B_1(0)}\frac{|h(x)|}{\rho_p(x)},\qquad \rho_p(x)=\sum_j\frac{\bar\varepsilon_{p,j}}{(\bar\varepsilon_{p,j}^2+|x-\xi_j|^2)^{3/2}}.\]
\end{lemma}
\begin{proof}
We use the potential expansion of \cite[Section 6.1, (6.4) and (6.13)]{Ianni-Luo-Yan-2025}.
For all parameters in the closures of these sets, write $\epsilon_j=\bar\varepsilon_{p,j}$. Unspecified integrals in $x$ and $y$ are over $B_1(0)$ and $\mathbb R^2$, respectively. For the signed sum, \eqref{sec5-projection} and the matching equations \eqref{Fpj} give
\[
\begin{gathered}
 \frac{s_jW_{\boldsymbol\alpha,p}(\xi_j+\epsilon_jy)}{a_{p,j}}
 =\alpha_j\left(p+U+\frac{w_0}p+\frac{w_1}{p^2}\right)
   +d_j+O(\epsilon_j|y|)+O(p^Ne^{-p/2}),\\
 \text{with}\quad |d_j|\le Cp^{-5/2},\qquad
 |\alpha_j^{p-1}-1|\le Cp^{-3/2}.
\end{gathered}
\]
Thus the calculation gives
\be\lab{sec5-potential-expansion}
    \epsilon_j^2p|W_{\boldsymbol\alpha,p}(\xi_j+\epsilon_j y)|^{p-1}=e^U\left(1+\frac{w_0-U-U^2/2}{p}\right)+\mathcal R_{p,j}(y), 
\ee
and
\be 
    \int_{B_{r/\epsilon_j}}(1+\log^2(2+|y|))|\mathcal R_{p,j}(y)|\,dy=o(p^{-1}).
\ee
To check the weighted remainder uniformly, fix $0<\delta<1$ and $L$ with $L(2-\delta)>1$. Taylor expansion on $|y|\le p^L$ and the profile decay on its complement yield, for fixed integers $K,N$,
\[
\begin{aligned}
 \int_{B_{p^L}}(1+\log(2+|y|))^2|\mathcal R_{p,j}|\,dy&\le C\big(p^{-3/2}+p^{-2}+p^Ne^{-p/4}\big),\\
 \int_{B_{r/\epsilon_j}\setminus B_{p^L}}(1+\log(2+|y|))^2|\mathcal R_{p,j}|\,dy&\le Cp^{-L(2-\delta)}(1+\log^Kp)=o(p^{-1}).
\end{aligned}
\]
Here the decay follows from the explicit profiles and \eqref{sec5-radial-expansion}.
For $Z(r)=(r^2-8)/(r^2+8)$ and $\chi=2+2U+rw_0'$, we have
    \[ \chi(r)=-8\log r+O(1),\qquad r\chi'(r)=-8+o(1).\]
The radial equations \eqref{3-12-0} and \eqref{eq-w1-correct} yield
    \[(-\Delta-e^U)Z=0,\qquad(-\Delta-e^U)\chi=-2e^U(w_0-U-U^2/2)Z.\]
Hence Green's identity gives
\[
 -2\int_{\mathbb R^2}e^U(w_0-U-U^2/2)Z^2\,dy
 =\lim_{R\to\infty}2\pi R(\chi Z'-Z\chi')(R)=16\pi,
\]
that is,
\be\lab{sec5-dilation-integral}
 \int_{\mathbb R^2}e^U(w_0-U-U^2/2)Z^2\,dy=-8\pi.
\ee
Moreover, direct computation gives
    \[\int e^UZ^2\,dy=\frac{8\pi}{3}.\]

We first obtain the multiplier estimate as in \cite[Lemma 6.2, (6.6)--(6.12) and (6.15)--(6.16)]{Ianni-Luo-Yan-2025}. If $v\in\mathcal O_e(\boldsymbol\xi)$ satisfies
\[-\Delta v-\sigma p|W_{\boldsymbol\alpha,p}|^{p-1}v=h+\sum_j\big(b_j(-\Delta Z_{p,j}) +c_j\partial_{\xi_{j,1}}(-\Delta Z_{p,j})\big),\]
with $0\le\sigma\le2$, then orthogonality  gives
    \[\|p|W_{\boldsymbol\alpha,p}|^{p-1}Z_{p,j}-p(-\Delta Z_{p,j})\|_1\le Cp^{-1},\]
    \[\left|\int\big(p|W_{\boldsymbol\alpha,p}|^{p-1}\partial_{\xi_{j,1}}Z_{p,j}-\partial_{\xi_{j,1}}(-\Delta Z_{p,j})\big)v\,dx\right|\le\frac{C\|v\|_\infty}{p^2\epsilon_j},\\\]
and
    \[ \sum_j\left(|b_j|+\frac{|c_j|}{\epsilon_j}\right)\le C\big(\|v\|_\infty+p\|h\|_*\big)+o(1)\sum_j\left(|b_j|+\frac{|c_j|}{\epsilon_j}\right).\]
Absorbing the last term proves \eqref{sec5-multiplier-bound}. 

To prove that \eqref{sec5-inverse} implies $\|v\|_\infty\le Cp\|h\|_*$, we now use the contradiction argument of \cite[Propositions 6.1 and 6.3]{Ianni-Luo-Yan-2025}. 
Suppose there is a sequence $\varphi_p\in\mathcal O_e(\boldsymbol\xi)$ with $\|\varphi_p\|_\infty=1$ such that
    \[-\Delta\varphi_p-\sigma_pp|W_{\boldsymbol\alpha,p}|^{p-1}\varphi_p=f_p+M_p,\]
with $0\le\sigma_p\le1+o(p^{-1})$ and $M_p=\sum_j\big(\beta_j(-\Delta Z_{p,j})+\gamma_j\partial_{\xi_{j,1}}(-\Delta Z_{p,j})\big)$, but $p\|f_p\|_*\to0$.
The multiplier bound \eqref{sec5-multiplier-bound} implies
\be\lab{sec5-multipliers}
 |\beta_j|\le C,\qquad |\gamma_j|\le C\epsilon_j,
 \qquad\|M_p\|_*\le Cp^{-1}.
\ee
Fix a sufficiently large $R$. On $\Omega_{p,R}=B_1(0)\setminus\bigcup_jB_{R\epsilon_j}(\xi_j)$, set
\[B_R(x)=\sum_j\left(2-\frac{R\epsilon_j}{|x-\xi_j|}\right),\]
which satisfies
\[2\le B_R\le4,\qquad-\Delta B_R=\sum_j\frac{R\epsilon_j}{|x-\xi_j|^3}.\]
The composite decay bound
\[
 p|W_{\boldsymbol\alpha,p}(x)|^{p-1}
 \le C\sum_k\epsilon_k^{-2}
 \left(1+\frac{|x-\xi_k|}{\epsilon_k}\right)^{-4+\delta},
 \qquad 0<\delta<1,
\]
gives
\[(-\Delta-\sigma_pp|W_{\boldsymbol\alpha,p}|^{p-1})B_R\ge cR\rho_p\quad\hbox{in }\Omega_{p,R},\]
and
\[1=\|\varphi_p\|_\infty\le C\max_j\|\varphi_p\|_{L^\infty(B_{R\epsilon_j}(\xi_j))}+\frac C R(p^{-1}+\|f_p\|_*).\]
Moreover, we have $\int\rho_p\,dx\le C$ and hence
    \[\|\varphi_p\|^2=\sigma_p\int p|W_{\boldsymbol\alpha,p}|^{p-1}\varphi_p^2\,dx+\int f_p\varphi_p\,dx\le C.\]
Then after passing to a subsequence, we get $\sigma_p\to\sigma\in[0,1]$ and
    \[\widehat\varphi_{p,j}(y):=\varphi_p(\xi_j+\epsilon_jy)\to
        \widehat\varphi_j\quad\hbox{in }C^1_{\mathrm{loc}}(\mathbb R^2),\]
where $\widehat\varphi_j\in\HR_e$ satisfies
    \[-\Delta\widehat\varphi_j=\sigma e^U\widehat\varphi_j,\qquad \int_{\mathbb R^2}|\nabla\widehat\varphi_j|^2<\infty,\]
and
    \[\max_j\|\widehat\varphi_j\|_\infty>0,\qquad\int e^U\widehat\varphi_j=0,\qquad \int e^U\partial_1U\,\widehat\varphi_j=0.\]
It follows from Lemmas \ref{lemma-3-10-1} and \ref{lemma-limit-kernel} that
    \[\sigma=1,\qquad \widehat\varphi_j=a_jZ.\]
Choose $j$ with $a_j\ne0$ and change the sign so that $a_j=a>0$.

Let $v$ be a radial solution of $(\Delta+e^U)v=e^UZ$. The radial ODE gives $v(r)=\frac43\log r+\kappa_v+o(1)$ as $r\to+\iy$.
For $j$ such that $a_j=a>0$, we define
    \[z_j(x)=Z\left(\frac{x-\xi_j}{\epsilon_j}\right),\qquad q_j(x)=\epsilon_j^{-2}e^{U((x-\xi_j)/\epsilon_j)},\]
and
    \[T_p=P\left[v\left(\frac{x-\xi_j}{\epsilon_j}\right)+\frac43\log\epsilon_j\,z_j
       +\left(\frac{8\pi}{3}H(\xi_j,\xi_j)+\kappa_v\right)(1-z_j)\right].\]
The expansion of projection $P$ yields
    \[T_p(\xi_j+\epsilon_jy)=\frac43\log\epsilon_j\,Z(y)+O(1),\]
and
    \[(-\Delta-q_j)T_p=-q_jz_j+q_j\frac{8\pi}{3}\big(H(x,\xi_j)-H(\xi_j,\xi_j)\big)+e_{p,j},\] 
with $\|e_{p,j}\|_1=o(1)$, so that $\|T_p\|_\infty\le Cp$.
We test the equation for $T_p$ with $\varphi_p$ and the equation for $\varphi_p$ with $T_p$. Subtracting the resulting identities gives
\be\lab{11-19-1}\begin{aligned}
 (1-\sigma_p)\int q_jT_p\varphi_p
 ={}&\sigma_p\int\big(p|W_{\boldsymbol\alpha,p}|^{p-1}-q_j\big)
          T_p\varphi_p+\int q_jz_j\varphi_p\\
 &+\int(f_p+M_p)T_p
   -\int q_j\frac{8\pi}{3}\big(H(x,\xi_j)-H(\xi_j,\xi_j)\big)\varphi_p
   +o(1).
\end{aligned}\ee
The definition of $Z_{p,k}$ in \eqref{Zpk-Wap} and the corresponding profile estimates give
    \[\left|\int(-\Delta Z_{p,k})T_p\right|\le Cp^{-1},\qquad
        \left|\int\partial_{\xi_{k,1}}(-\Delta Z_{p,k})\,T_p\right|\le C,\quad\forall~k=1,2.\]
By the above estimates of $W_{\boldsymbol\alpha,p},\varphi_p$ and $T_p$, we get
    \[\int_{B_r(\xi_{j'})}p|W_{\boldsymbol\alpha,p}|^{p-1}T_p\varphi_p=o(1),\quad j'\ne j,\]
    \[\int q_j\frac{8\pi}{3}\big(H(x,\xi_j)-H(\xi_j,\xi_j)\big)\varphi_p=o(1),\]
and
    \[\left|\int(f_p+M_p)T_p\right|\le Cp\|f_p\|_*+\frac C p\sum_k|\beta_k|
       +C\sum_k|\gamma_k|=o(1).\]
By \eqref{sec5-potential-expansion}, \eqref{sec5-dilation-integral}, and $p^{-1}\log\epsilon_j\to-1/4$, we get
    \[\int q_jz_j\varphi_p=\frac{8\pi a}{3}+o(1),\]
    \[\int q_jT_p\varphi_p =\frac43\log\epsilon_j\left(\frac{8\pi a}{3}+o(1)\right)+O(1),\]
and
    \[ \sigma_p\int\big(p|W_{\boldsymbol\alpha,p}|^{p-1}-q_j\big)T_p\varphi_p
        =\frac{4\sigma_pa\log\epsilon_j}{3p}\int e^U(w_0-U-U^2/2)Z^2+o(1) =\frac{8\pi a}{3}+o(1).\]
Consequently, we get from \eqref{11-19-1} that
    \[\frac{16\pi a}{3}=\frac43(1-\sigma_p)\log\epsilon_j\left(\frac{8\pi a}{3}+o(1)\right)+o(1)\le o(1),\]
which is a contradiction since we assume $a_j=a>0$. 

If \eqref{sec5-coercive} failed, minimizers of
\[
 \zeta_p=\inf_{\substack{v\in\mathcal O_e(\boldsymbol\xi)\\
                         \int p|W_{\boldsymbol\alpha,p}|^{p-1}v^2>0}}
       \frac{\|v\|^2}{\int p|W_{\boldsymbol\alpha,p}|^{p-1}v^2}
\]
would give the excluded sequence $\vp_p$ with $\sigma_p=\zeta_p\le1+o(p^{-1})$ and $f_p=0$.
If the uniform $L^\infty$ estimate failed for solutions $v_p$ of \eqref{sec5-inverse}, normalizing by $\|v_p\|_\infty$ would give the excluded sequence with
\[ 
\sigma_p=t_p\in[0,1],\qquad
\varphi_p=\frac{v_p}{\|v_p\|_\infty},\qquad
f_p=\frac{h_p}{\|v_p\|_\infty},\qquad
p\|f_p\|_*=\frac{p\|h_p\|_*}{\|v_p\|_\infty}\longrightarrow0.
\]
\end{proof}

By Lemma \ref{sec5-linear-lemma}, the operator $K_p=\Pi_{\boldsymbol\xi}(-\Delta)^{-1}(p|W_{\boldsymbol\alpha,p}|^{p-1}\,\cdot)$ is compact on $\mathcal O_e(\boldsymbol\xi)$ and satisfies
\[
\langle(I-tK_p)v,v\rangle\ge\left(1-t+\frac{ct}{p}\right)\|v\|^2,\qquad \operatorname{ind}(I-tK_p)=0\quad(0\le t\le1).
\]
Thus $I-tK_p$ is invertible.

Let $\mathcal B_p(\boldsymbol\alpha,\boldsymbol\xi,w)$ have the first four components of $\mathcal A_p$ evaluated at $w=0$, and third block equal to $w$. Thus
\[\begin{aligned}
 \mathcal B_{p,1,j}&=\langle J_p'(W_{\boldsymbol\alpha,p}),Z_{p,j}\rangle,\\
 \mathcal B_{p,2,j}&=\mathcal P_j(W_{\boldsymbol\alpha,p}),\\
 \mathcal B_{p,3}&=w,
\end{aligned}\]
and $\BR_p=(\BR_{p,1,1},\BR_{p,1,2},\BR_{p,2,1},\BR_{p,2,2},\BR_{p,3})$.
At $\boldsymbol\alpha=\boldsymbol1$, put $V=G(\xi_1,\cdot)-G(\xi_2,\cdot)$.
This function is harmonic on the annulus supporting $\nabla\eta_j$. Integration by parts in \eqref{sec5-stress} gives one half of the boundary expression $Q_j(V,V)$, where $Q_j$ is given by \eqref{Qj-1}. The formulas for $Q_j$ give
\[
 Q_j(V,V)=\partial_1R(\xi_j)-2\partial_{\xi_{j,1}}G(\xi_1,\xi_2)
        =-\partial_{\xi_{j,1}}\Psi(\boldsymbol\xi).
\]
Using \eqref{sec5-outer-W}, we obtain, uniformly on the fixed parameter neighborhoods,
\be\lab{sec5-B2}
 \mathcal B_{p,2,j}
 =-\frac{32\pi^2e}{p^2}\partial_{\xi_{j,1}}\Psi(\boldsymbol\xi)
 +o(p^{-2}).
\ee
The same formula holds for $|\alpha_j-1|\le p^{-5/2}$.

\begin{lemma}\lab{sec5-residual-lemma}
The residual of $W_{\boldsymbol\alpha,p}$ can be decomposed as
\be\lab{sec5-residual}
 -\Delta W_{\boldsymbol\alpha,p}-|W_{\boldsymbol\alpha,p}|^{p-1}W_{\boldsymbol\alpha,p}
 =\sum_j b_j(-\Delta Z_{p,j})+h_{\boldsymbol\alpha,p},
 \qquad \|h_{\boldsymbol\alpha,p}\|_*\le Cp^{-4}.
\ee
The amplitude components $\mathcal B_{p,1,j}$ also satisfy
\be\lab{sec5-B1}
 \mathcal B_{p,1,j}=8\pi p a_{p,j}^2(\alpha_j-\alpha_j^p)+O\left(p^{-2}\sum_i|\alpha_i-\alpha_i^p|+p^{-4}\right).
\ee
\end{lemma}
\begin{proof}
Near $\xi_j$, we have
\[
\begin{aligned}
 s_jW_{\boldsymbol\alpha,p}
 &=a_{p,j}\left[\alpha_j(p+U+w_0/p+w_1/p^2)+d_j\right]+O(a_{p,j}|x-\xi_j|)+O(p^Ne^{-p/2}),
\end{aligned}
\]
with $|d_j|\le C\max_i|\alpha_i-1|$.
Set $\widehat\alpha_j=\alpha_j+d_j/p$ and  $b_j=\alpha_j-\widehat\alpha_j^p$.
In a fixed small ball centered at $\xi_j$, the above expansion and the equations for the corrections $\sum_j b_j(-\Delta Z_{p,j})$ show that the remaining part is bounded by
\[
 \frac{Ca_{p,j}}{\bar\varepsilon_{p,j}^2}e^U(1+|y|)^\delta\left[p^{-3}+p^{-1}\max_i|\alpha_i-1|\right].
\]
For $0<\delta<1$, division by $\rho_p$ bounds the first term by $Cp^{-4}$ and the second by $Cp^{-9/2}$. The region outside the fixed balls is treated by \eqref{sec5-outer-W} and \eqref{sec5-radial-expansion}. This proves \eqref{sec5-residual}.

For \eqref{sec5-B1}, pair the original residual with $Z_{p,j}$. In the $j$-th ball its leading coefficient is $\alpha_j-\alpha_j^p$, and the local integral equals $8\pi p a_{p,j}^2(1+O(p^{-1}))$. Contributions from the other peak ($j'\neq j$) have an additional factor $p^{-1}$, since $Z_{p,j}=O(p^{-1})$ there. The shift $d_j$ and the profile correction errors give the displayed remainder.
\end{proof}

\begin{proof}[Proof of Proposition \ref{thm-degree}]
We use the homotopy of \cite[Proposition~7.5]{Ianni-Luo-Yan-2025} on a fixed space. Define the open sets $\widehat{\mathcal S}_p^\pm$ by taking the bounds in $\mathcal S_p^\pm$ to be strict and letting $w$ range over $\mathcal X_e$. Write $F_p=I-\mathcal T_p$ and let $\Pi_{\boldsymbol\xi}$ denote the Dirichlet orthogonal projection onto $\mathcal O_e(\boldsymbol\xi)$. Extend $\mathcal A_p$ to $\widehat{\mathcal A}_p$ by taking its last block to be
\[
 \Pi_{\boldsymbol\xi}F_p(W_{\boldsymbol\alpha,p}+w)+(I-\Pi_{\boldsymbol\xi})w
 =w-\Pi_{\boldsymbol\xi}\mathcal T_p(W_{\boldsymbol\alpha,p}+w).
\]
With $\mathcal B_p$ given by the same formulas on the enlarged sets, the homotopy $t\widehat{\mathcal A}_p+(1-t)\mathcal B_p$ consists of compact perturbations of the identity. At every homotopy zero, its last block forces $(I-\Pi_{\boldsymbol\xi})w=0$. Thus we may use the following boundary estimates for $\mathcal A_{p,t}=t\mathcal A_p+(1-t)\mathcal B_p$ with $w\in\mathcal O_e(\boldsymbol\xi)$.
At such a zero, the last block of $\mathcal A_{p,t}$ gives the equation
\be\lab{11-19-2}
 -\Delta w-tp|W_{\boldsymbol\alpha,p}|^{p-1}w=-t h_{\boldsymbol\alpha,p}+N_t(w)
 +\sum_j\left[\beta_j(-\Delta Z_{p,j})+
 \gamma_j(-\Delta\partial_{\xi_{j,1}}Z_{p,j})\right],
\ee
where $h_{\boldsymbol\alpha,p}$ is given by \eqref{sec5-residual}. Taylor's formula for $f(u)=|u|^{p-1}u$, with $p$ large, gives
\[
 \|N_t(w)\|_*\le Cp\|w\|_\infty^2.
\]
Indeed, the estimate of $N_t(w)$ near a peak follows from
$p(p-1)(|W_{\boldsymbol\alpha,p}|+|w|)^{p-2}\le Cp\rho_p$,
and off the fixed peak balls $N_t(w)$ is smaller than every negative power of $p$. Since $p^2\|w\|_\infty\le p^{-\tau_0}$, \eqref{sec5-inverse} implies
\[
 \|w\|_\infty\le Cp^{-3}+Cp^2\|w\|_\infty^2
 \le Cp^{-3}+Cp^{-\tau_0}\|w\|_\infty.
\]
This gives $\|w\|_\infty\le Cp^{-3}$ for $p$ large enough.
Since $w\in\OR_e$, by \eqref{sec5-coercive} and \eqref{11-19-2} we get
\[
 \frac c p\|w\|^2
 \le\left|\int h_{\boldsymbol\alpha,p}w\,dx\right|
     +\left|\int N_t(w)w\,dx\right|
 \le Cp^{-4}\|w\|_\infty+Cp\|w\|_\infty^3.
\]
It follows that $\|w\|\le Cp^{-3}$. The multiplier and exterior estimates follow the argument of \cite[Proposition 7.5, Step 1]{Ianni-Luo-Yan-2025}. Indeed, for $f_t=-t h_{\boldsymbol\alpha,p}+N_t(w)$, we have $\|f_t\|_*\le Cp^{-4}$.
Then by \eqref{sec5-multiplier-bound} we get
\[
 \sum_j\left(|\beta_j|+\frac{|\gamma_j|}{\bar\varepsilon_{p,j}}\right)\le C\big(\|w\|_\infty+p\|f_t\|_*\big)\le Cp^{-3}.
\]
Let $\Omega_r=B_1(0)\setminus\bigcup_jB_r(\xi_j)$. By \eqref{sec5-outer-W} and \eqref{11-19-2} , we have  $\|\Delta w\|_{L^q(\Omega_{r/2})}\le Cp^{-3}$ for fixed $q>2$. 
Then elliptic estimates yield
\[
 \|w\|_{C^1(\overline{\Omega_r})}
 \le C\big(\|w\|_\infty+\|\Delta w\|_{L^q(\Omega_{r/2})}\big)
 \le Cp^{-3}.
\]
By orthogonality we get
\[
 \mathcal A_{p,1,j}-\mathcal B_{p,1,j}
 =-\int p|W_{\boldsymbol\alpha,p}|^{p-1}Z_{p,j}w\,dx
   +O(p\|w\|_\infty^2)=O(p^{-3}).
\]
Then, the equation $\mathcal A_{p,t,1,j}=0$, together with \eqref{sec5-B1}, yields
\[
 \sum_j|\alpha_j-\alpha_j^p|\le Cp^{-2},\qquad
 |\alpha_j-1|\le Cp^{-3}.
\]
Thus no amplitude or remainder boundary can contain a homotopy zero.

On the support of $\nabla\eta_j$, the gradients of $W_{\boldsymbol\alpha,p}$ and $w$ have orders $p^{-1}$ and $p^{-3}$, respectively. Formula \eqref{sec5-stress} therefore gives
\[
 \mathcal A_{p,2,j}-\mathcal B_{p,2,j}=O(p^{-4}).
\]
By \eqref{sec5-B2}, a homotopy zero satisfies $\nabla_{\boldsymbol\xi}\Psi(\boldsymbol\xi)=o(1)$, with derivatives restricted to the two axial coordinates. The Hessian restricted to two axial coordinates at $\pm\boldsymbol X$ is positive definite. Taking $d_0$ small enough, these are the only restricted critical points in the two closed boxes, and their boundaries satisfy $|\nabla_{\boldsymbol\xi}\Psi|\ge c(d_0)>0$. Hence the location boundary is also excluded.

Lemma \ref{lemma-3-12-1} identifies all zeros of $\mathrm{id}-\mathcal T_p$ with those of $\widehat{\mathcal A}_p$.
We compare the local degrees using the tests of \cite[Proposition~7.1 and Remark~7.2]{Ianni-Luo-Yan-2025}. At a zero $(q,w)$, write $q=(\alpha_1,\alpha_2,\xi_{1,1},\xi_{2,1})$, $W=W_{\boldsymbol\alpha,p}$, $\Pi=\Pi_{\boldsymbol\xi}$, $u=W+w$, and $L=DF_p(u)$. Use the ordered families $B=(Z_{p,1},Z_{p,2},Y_{p,1},Y_{p,2})$ and $T=(Z_{p,1},Z_{p,2},\eta_1\partial_1u,\eta_2\partial_1u)$ from Lemma~\ref{lemma-3-12-1}. By the stress identity and density, $D_u\mathcal P_j(u)[v]=\langle Lv,\eta_j\partial_1u\rangle$, and the derivative in the cutoff vanishes since $F_p(u)=0$. Write $\mathcal X_e=O\oplus E$ with $O=\Pi\mathcal X_e$ and $E=\operatorname{span}B$. Precomposing $D\widehat{\mathcal A}_p$ with the shear
\[
 (\dot q,h,r)\longmapsto(\dot q,h+r+D\Pi[\dot q]w),
 \qquad h\in O,\quad r\in E .
\]
gives an upper triangular operator with diagonal blocks $RLC$ and $I_E$. Here $D\Pi[\dot q]w\in E$ follows from $\Pi^2=\Pi$, and
\[
 C(\dot q,h)=DW[\dot q]+D\Pi[\dot q]w+h,\qquad
 Rv=\big((\langle v,T_i\rangle)_{i=1}^4,\Pi v\big).
\]
Identify $\mathbb R^4\times O$ with $\mathcal X_e$ by $U_B(b,h)=\sum_i b_iB_i+h$. In $RLC=(RU_B)(U_B^{-1}LU_B)(U_B^{-1}C)$, the outer factors are block triangular with the identity on $O$. The finite block of $RU_B$ has positive determinant by the normalized test matrix in Lemma~\ref{lemma-3-12-1}, whose limit is $\operatorname{diag}(1,1,-1,-1)$. For $U_B^{-1}C$, differentiating $\langle\Pi w,B_i\rangle=0$ with $w$ fixed identifies the finite block with the inverse Gram matrix times $-D\Phi_p$, up to positive row and column scalings. Its determinant is positive by the Jacobian estimate in the same lemma. The outer factors and the shear thus deform to the identity through invertible compact perturbations of the identity. The off-diagonal block above can be scaled to zero, and conjugation preserves the index of $L$, which is invertible by Proposition~\ref{thm-index} and compactness. Hence the two local degrees coincide.

Thus, using \eqref{sec5-degree-spaces} and additivity and homotopy invariance, we get
\[
 \deg(\mathrm{id}-\mathcal T_p,\DR_p\cap\mathcal H_e,0)
 =\deg(\widehat{\mathcal A}_p,\widehat{\mathcal S}_p^+\cup\widehat{\mathcal S}_p^-,0)
 =\deg(\mathcal B_p,\widehat{\mathcal S}_p^+\cup\widehat{\mathcal S}_p^-,0).
\]
The estimates \eqref{sec5-B1} and \eqref{sec5-B2} permit a further boundary-nonvanishing homotopy to
\[
 (\boldsymbol\alpha,\boldsymbol\xi,w)\longmapsto
 \left(\frac{8\pi e}{p}(\alpha_1-\alpha_1^p,\alpha_2-\alpha_2^p),
 -\frac{32\pi^2e}{p^2}\nabla_{\boldsymbol\xi}\Psi(\boldsymbol\xi),w\right).
\]
The amplitude block has degree $+1$, since both derivatives at $\boldsymbol\alpha=(1,1)$ are negative. Each location block has degree $+1$, since its derivative is a negative scalar multiple of the positive definite two-dimensional restricted Hessian. The identity on the orthogonal complement has degree $+1$. The two disjoint location boxes therefore contribute $1+1=2$.
\end{proof}

\vskip0.2in
\section{Antisymmetry, uniqueness and nondegeneracy}\lab{sec-6}

\begin{proof}[Proof of Theorem \ref{thm-symmetry}]
Lemma~\ref{lemma-3-6-2} and Proposition~\ref{lemma-1} give unique global extrema and symmetry about the line through the extrema and the origin. Rotate this line onto the $x_1$-axis. Then $u_p\in\ER_p\cap\CR_p\cap\HR_e$.

By Proposition~\ref{thm-index}, every solution in this set is nondegenerate in $\HR_e$ and has Morse index two there. 
The solution set in $\DR_p\cap\HR_e$ is compact because $\DR_p$ is bounded, $\TR_p$ is compact and there are no boundary fixed points. Nondegeneracy makes each solution isolated, so there are finitely many solutions.
The existence result of Esposito, Musso and Pistoia, together with Lemma~\ref{lem-mass-energy} and Proposition~\ref{lemma-1}, gives at least one such solution. Let $k_p$ denote the number of these solutions.
The local index formula and additivity of the degree give
\[
\deg(I-\TR_p,\DR_p\cap\HR_e,0)
=\sum_{u\in\ER_p\cap\CR_p\cap\HR_e}(-1)^{m_e(u)}=k_p.
\]
Proposition~\ref{thm-degree} therefore implies
\be\lab{3-12-20-1}
\sharp(\ER_p\cap\CR_p\cap\HR_e)=k_p=2.
\ee
Both $-u_p$ and $u_p(-\cdot)$ belong to this set. The first differs from $u_p$ because $u_p$ is nonzero. The second also differs from $u_p$, since otherwise the unique positive maximum point would equal its antipode and hence would be the origin, contrary to its convergence to a nonzero concentration point. Thus
\be\lab{odd}
-u_p(x_1,x_2)=u_p(-x_1,-x_2)=u_p(-x_1,x_2),
\qquad (x_1,x_2)\in B_1(0),
\ee
where the last equality uses evenness in $x_2$. This is antisymmetry about the $x_2$-axis.
Rotate $u_p$ by $\pi$ if necessary so that its maximum point lies on the positive $x_1$-axis. By Corollary~\ref{cor-angular-monotonicity}, $u_p(r,\theta)$ is strictly decreasing on $0<\theta<\pi$ and vanishes at $\theta=\pi/2$ by \eqref{odd}. It is therefore positive for $x_1>0$ and negative for $x_1<0$ in the upper half-disk. Evenness in $x_2$ and continuity give the same signs throughout the respective open half-disks, so the interior nodal set is $B_1(0)\cap\{x_1=0\}$.
\end{proof}

\begin{proof}[Proof of Theorem \ref{thm-uniqueness}]
Every solution in $\ER_p\cap\CR_p$ can be rotated into $\HR_e$ by Proposition~\ref{lemma-1}. Equations~\eqref{3-12-20-1} and \eqref{odd} identify the two solutions in that subspace as $u_p$ and $u_p(-\cdot)$, which differ by a rotation through $\pi$. Hence there is exactly one rotation class.
Choose a representative whose maximum point lies on the positive $x_1$-axis.
To prove the assertions about the full linearized operator, decompose
\[H_0^1(B_1(0))=\HR_e\oplus\HR_o,\qquad
      \HR_o=\{v\in H_0^1(B_1(0)):v(x_1,-x_2)=-v(x_1,x_2)\}.\]
The two subspaces are orthogonal for both the Dirichlet inner product and the bilinear form
\[\QR_p(v,w)=\int_{B_1(0)}\nabla v\cdot\nabla w\rd x
      -p\int_{B_1(0)}|u_p|^{p-1}vw\rd x,\]
because $|u_p|^{p-1}$ is even in $x_2$. Proposition \ref{thm-index} gives two negative directions and a trivial kernel on $\HR_e$.

Following the use of the angular derivative on a half-disk in \cite[proof of Theorem~1.3]{Aftalion-Pacella-2004}, let $B_1^+=B_1(0)\cap\{x_2>0\}$ and set $z_p=-\partial_\theta u_p$. By Corollary~\ref{cor-angular-monotonicity}, $z_p>0$ in $B_1^+$. Rotational invariance and the symmetry of $u_p$ give
\[-\Delta z_p=p|u_p|^{p-1}z_p\quad\text{in }B_1^+,
      \qquad z_p=0\quad\text{on }\partial B_1^+.\]
Let $\mu_1$ be the first eigenvalue of the weighted Dirichlet problem on $B_1^+$, and let $\varphi_1>0$ be a corresponding eigenfunction. Testing the equations for $z_p$ and $\varphi_1$ against each other yields
\[(1-\mu_1)p\int_{B_1^+}|u_p|^{p-1}z_p\varphi_1\rd x=0.\]
The integral is strictly positive, so $\mu_1=1$. The first eigenvalue is simple on the connected domain $B_1^+$. Restriction to $B_1^+$ identifies $\HR_o$ with $H_0^1(B_1^+)$ by odd reflection. Consequently, $\QR_p$ is nonnegative on $\HR_o$, and its kernel there is spanned by $\partial_\theta u_p$.
Combining the two invariant subspaces, we obtain
\[m(u_p)=2,\qquad
      \ker\bigl(-\Delta-p|u_p|^{p-1}\bigr)=\R\partial_\theta u_p.\qedhere\]
\end{proof}

\begin{proof}[Proof of Corollary \ref{cor-1}]
Every family of least-energy nodal solutions is a low-energy family by \cite{Grossi-Grumiau-Pacella-2013}. Lemma~\ref{lem-mass-energy} therefore places it in $\ER_p$ for all sufficiently large $p$. Theorem~\ref{thm-uniqueness} gives uniqueness up to rotations.
\end{proof}

\begin{proof}[Proof of Corollary \ref{cor-2}]
Let $\{u_p\}$ be a low-energy nodal family, and choose a family $\{v_p\}$ of least-energy nodal solutions. By Lemma~\ref{lem-mass-energy}, both lie in $\ER_p\cap\CR_p$ for all sufficiently large $p$. Theorem~\ref{thm-uniqueness} identifies $u_p$ and $v_p$ up to rotations. In particular, this applies to the family constructed in Theorem~\ref{Esposito-Musso-Pistoia}.
\end{proof}

\vs
\noindent {\bf Acknowledgements}
The research of Z. Chen is supported by the National Key \text{R\&D} Program of China (Grant 2023YFA1010002) and the National Natural Science Foundation of China (No. 1222109).
The research of H. Li is supported by the Beijing Institute of Technology Research Fund Program for Young Scholars and the National Natural Science Foundation of China (No. 12601182).


\begin{thebibliography}{10}

\bibitem{Adimurthi-Grossi-2004}
\textsc{Adimurthi and M.~Grossi},
\emph{Asymptotic estimates for a two-dimensional
  problem with polynomial nonlinearity}, Proc. Amer. Math. Soc., 132 (2004),
  pp.~1013--1019.

\bibitem{Aftalion-Pacella-2004}
\textsc{A.~Aftalion and F.~Pacella},
\emph{Qualitative properties of nodal solutions of semilinear elliptic equations
  in radially symmetric domains}, C. R. Math. Acad. Sci. Paris, 339 (2004),
  pp.~339--344.

\bibitem{Bartsch-Weth-2003}
\textsc{T.~Bartsch and T.~Weth},
\emph{A note on additional properties of sign changing solutions to superlinear
  elliptic equations}, Topol. Methods Nonlinear Anal., 22 (2003), pp.~1--14.

\bibitem{Bartsch-Weth-Willem-2005}
\textsc{T.~Bartsch, T.~Weth, and M.~Willem},
\emph{Partial symmetry of least energy
  nodal solutions to some variational problems}, J. Anal. Math., 96 (2005),
  pp.~1--18.

\bibitem{Bonheure-Bouchez-Grumiau-VanSchaftingen-2008}
\textsc{D.~Bonheure, V.~Bouchez, C.~Grumiau, and J.~Van Schaftingen},
\emph{Asymptotics and symmetries of least energy nodal solutions of Lane--Emden
  problems with slow growth}, Commun. Contemp. Math., 10 (2008), pp.~609--631.

\bibitem{Cao-Guo-Peng-Yan-2019}
\textsc{D.~Cao, Y.~Guo, S.~Peng, and S.~Yan},
\emph{Local uniqueness for vortex
  patch problem in incompressible planar steady flow}, J. Math. Pures Appl.
  (9), 131 (2019), pp.~251--289.

\bibitem{Chen-Lin-2003}
\textsc{C.-C. Chen and C.-S. Lin},
\emph{Topological degree for a mean field
  equation on {R}iemann surfaces}, Comm. Pure Appl. Math., 56 (2003),
  pp.~1667--1727.

\bibitem{Chen-Cheng-Zhao-2023}
\textsc{Z.~Chen, Z.~Cheng, and H.~Zhao},
\emph{Asymptotic behavior of least energy nodal solutions for biharmonic
  Lane--Emden problems in dimension four}, Calc. Var. Partial Differential
  Equations, 62 (2023), Paper No.~205.

\bibitem{Chen-Cheng-Zhao-2024}
\textsc{Z.~Chen, Z.~Cheng, and H.~Zhao},
\emph{Sharp estimates and non-degeneracy
  of low energy nodal solutions for the {L}ane-{E}mden equation in dimension
  two}, Calc. Var. Partial Differential Equations, 63 (2024), Paper No.~234, 28 pp.

\bibitem{DeMarchis-Grossi-Ianni-Pacella-2018}
\textsc{F.~De~Marchis, M.~Grossi, I.~Ianni, and F.~Pacella},
\emph{{$L^\infty$}-norm and energy quantization for the planar {L}ane--{E}mden
  problem with large exponent}, Arch. Math. (Basel), 111 (2018), pp.~421--429.

\bibitem{DeMarchis-Grossi-Ianni-Pacella-2019}
\textsc{F.~De~Marchis, M.~Grossi, I.~Ianni, and F.~Pacella},
\emph{Morse index and
  uniqueness of positive solutions of the {L}ane-{E}mden problem in planar
  domains}, J. Math. Pures Appl. (9), 128 (2019), pp.~339--378.

\bibitem{DeMarchis-Ianni-Pacella-2015}
\textsc{F.~De~Marchis, I.~Ianni, and F.~Pacella},
\emph{Asymptotic analysis and
  sign-changing bubble towers for {L}ane-{E}mden problems}, J. Eur. Math. Soc.
  (JEMS), 17 (2015), pp.~2037--2068.

\bibitem{DeMarchis-Ianni-Pacella-2016}
\textsc{F.~De~Marchis, I.~Ianni, and F.~Pacella},
\emph{Morse index and
  sign-changing bubble towers for {L}ane-{E}mden problems}, Ann. Mat. Pura
  Appl. (4), 195 (2016), pp.~357--369.

\bibitem{DeMarchis-Ianni-Pacella-2017-2}
\textsc{F.~De~Marchis, I.~Ianni, and F.~Pacella},
\emph{Asymptotic profile
  of positive solutions of {L}ane--{E}mden problems in dimension two}, J. Fixed
  Point Theory Appl., 19 (2017), pp.~889--916.

\bibitem{DeMarchis-Ianni-Pacella-2017}
\textsc{F.~De~Marchis, I.~Ianni, and F.~Pacella},
\emph{Exact {M}orse index
  computation for nodal radial solutions of {L}ane-{E}mden problems}, Math.
  Ann., 367 (2017), pp.~185--227.

\bibitem{Esposito-Musso-Pistoia-2007}
\textsc{P.~Esposito, M.~Musso, and A.~Pistoia},
\emph{On the existence and profile
  of nodal solutions for a two-dimensional elliptic problem with large exponent
  in nonlinearity}, Proc. Lond. Math. Soc. (3), 94 (2007), pp.~497--519.

\bibitem{Fazekas-Pacella-Plum-2022}
\textsc{B.~Fazekas, F.~Pacella, and M.~Plum},
\emph{Approximate nonradial solutions for the Lane--Emden problem in the ball},
  Adv. Nonlinear Anal., 11 (2022), pp.~268--284.

\bibitem{Gidas-Ni-Nirenberg-1979}
\textsc{B.~Gidas, W.-M. Ni, and L.~Nirenberg},
\emph{Symmetry and related properties via the maximum principle},
  Comm. Math. Phys., 68 (1979), pp.~209--243.

\bibitem{Gladiali-Stegel-2022}
\textsc{F.~Gladiali and G.~Stegel},
\emph{A note on nonradial nodal solutions to
  the {H}\'enon problem in the disc}, Indiana Univ. Math. J., 71 (2022),
  pp.~715--733.

\bibitem{Grossi-Grumiau-Pacella-2013}
\textsc{M.~Grossi, C.~Grumiau, and F.~Pacella},
\emph{Lane-{E}mden problems:
  asymptotic behavior of low energy nodal solutions}, Ann. Inst. H.
  Poincar\'e{} C Anal. Non Lin\'eaire, 30 (2013), pp.~121--140.

\bibitem{Grossi-Grumiau-Pacella-2014}
\textsc{M.~Grossi, C.~Grumiau, and F.~Pacella},
\emph{Lane {E}mden
  problems with large exponents and singular {L}iouville equations}, J. Math.
  Pures Appl. (9), 101 (2014), pp.~735--754.

\bibitem{Grossi-Ianni-Luo-Yan-2022}
\textsc{M.~Grossi, I.~Ianni, P.~Luo, and S.~Yan},
\emph{Non-degeneracy and local
  uniqueness of positive solutions to the {L}ane-{E}mden problem in dimension
  two}, J. Math. Pures Appl. (9), 157 (2022), pp.~145--210.

\bibitem{Grumiau-Troestler-2010}
\textsc{C.~Grumiau and C.~Troestler},
\emph{Oddness of least energy nodal solutions
  on radial domains}, Electron. J. Differ. Equ. Conf., 18 (2010), pp.~23--31.

\bibitem{Ianni-Luo-Yan-2025}
\textsc{I.~Ianni, P.~Luo, and S.~Yan},
\emph{Morse index, topological degree and
  local uniqueness of multi-spikes solutions to the Lane-Emden problem in
  dimension two}, arXiv:2506.12905, 2025.

\bibitem{Ianni-Saldana-2021}
\textsc{I.~Ianni and A.~Salda\~na},
\emph{Sharp asymptotic behavior of radial
  solutions of some planar semilinear elliptic problems}, J. Differential
  Equations, 304 (2021), pp.~102--164.

\bibitem{Kamburov-Sirakov-2018}
\textsc{N.~Kamburov and B.~Sirakov},
\emph{Uniform a priori estimates for positive
  solutions of the {L}ane--{E}mden equation in the plane}, Calc. Var. Partial
  Differential Equations, 57 (2018), Paper No.~164, 8 pp.

\bibitem{Li-1999}
\textsc{Y.~Y. Li},
\emph{Harnack type inequality: the method of moving planes},
  Comm. Math. Phys., 200 (1999), pp.~421--444.

\bibitem{Ni-Takagi-1991}
\textsc{W.-M. Ni and I.~Takagi},
\emph{On the shape of least-energy solutions to a
  semilinear {N}eumann problem}, Comm. Pure Appl. Math., 44 (1991),
  pp.~819--851.

\bibitem{Pacella-Weth-2007}
\textsc{F.~Pacella and T.~Weth},
\emph{Symmetry of solutions to semilinear elliptic equations via Morse index},
  Proc. Amer. Math. Soc., 135 (2007), pp.~1753--1762.

\bibitem{Ren-Wei-94}
\textsc{X.~Ren and J.~Wei},
\emph{On a two-dimensional elliptic problem with large
  exponent in nonlinearity}, Trans. Amer. Math. Soc., 343 (1994), pp.~749--763.

\bibitem{Ren-Wei-96}
\textsc{X.~Ren and J.~Wei},
\emph{Single-point
  condensation and least-energy solutions}, Proc. Amer. Math. Soc., 124 (1996),
  pp.~111--120.

\bibitem{Thizy-2019}
\textsc{P.-D. Thizy},
\emph{Sharp quantization for {L}ane--{E}mden problems in
  dimension two}, Pacific J. Math., 300 (2019), pp.~491--497.

\bibitem{Wei-Winter-2005}
\textsc{J.~Wei and M.~Winter},
\emph{Symmetry of nodal solutions for singularly
  perturbed elliptic problems on a ball}, Indiana Univ. Math. J., 54 (2005),
  pp.~707--741.
\end{thebibliography}
\end{document}